\documentclass[12pt]{article}
\usepackage{amsmath}
\usepackage{amssymb}
\numberwithin{equation}{section}

\title{ Integrals over classical Groups,
 Random permutations, Toda and Toeplitz lattices }

\author{
M. Adler\thanks{ Department of Mathematics, Brandeis
University, Waltham, Mass 02454, USA. E-mail:
adler@math.brandeis.edu.  The support of a National
Science Foundation grant \# DMS-98-4-50790 is
gratefully acknowledged.}~~~~~~ P. van
Moerbeke\thanks{ Department of Mathematics,
Universit\'e de Louvain, 1348 Louvain-la-Neuve,
Belgium and Brandeis University, Waltham, Mass 02454,
USA. E-mail: vanmoerbeke@geom.ucl.ac.be and
@math.brandeis.edu. The support of a National Science
Foundation grant \# DMS-98-4-50790, a Nato, a FNRS and
a Francqui Foundation grant is gratefully
acknowledged.}}

\date{November 30, 1999}

\catcode`\@=11
\let\c@equation=\relax
\newcounter{equation}[subsection]

\catcode`\@=12

\newcommand{\MAT}[1]{\left(\begin{array}{*#1c}}
\newcommand{\mat}{\end{array}\right)}
\newcommand{\qed}{\leavevmode\unskip\nobreak\penalty200\hskip2pt\null
\nobreak\hfill\rule{1.1ex}{1.1ex}%\parfillskip=0pt
\medbreak
}

\newcommand{\pp}{\ldots}

\newcommand{\JR}{{\cal J}}
\newcommand{\LR}{{\cal L}}

\newcommand{\VR}{{\cal V}}

\newcommand{\BC}{{\mathbb C}}

\newcommand{\BX}{{\mathbb X}}

\newcommand{\BZ}{{\mathbb Z}}

\newcommand{\iy}{\infty}
\newcommand{\pl}{\partial}
\newcommand{\al}{\alpha}

\newenvironment{proof}{\medskip\noindent{\it Proof:\/} }{\qed}
\newenvironment
        {remark}{\medskip\noindent\underline{\it Remark:\/} }{\medbreak}

\newcommand{\la}{\langle}
\newcommand{\ra}{\rangle}

\newcommand{\dt}{\delta}
\newcommand{\Dt}{\Delta}

\newcommand{\BR}{{\mathbb R}}
\newcommand{\lb}{\lambda}
\newcommand{\Lb}{\Lambda}

\newcommand{\BJ}{{\mathbb J}}

\newcommand{\diag}{\operatorname{diag}}

\def\be#1\ee{\begin{equation}#1\end{equation}}
\def\bea#1\eea{\begin{eqnarray}#1\end{eqnarray}}
\def\bean#1\eean{\begin{eqnarray*}#1\end{eqnarray*}}

\newcommand{\Tr}{\operatorname{\rm Tr}}

\newtheorem{definition}{Definition}[section]
\newtheorem{theorem}[definition]{Theorem}
\newtheorem{lemma}[definition]{Lemma}
\newtheorem{corollary}[definition]{Corollary}
\newtheorem{proposition}[definition]{Proposition}

\makeatletter
\def\ps@X{\let\@mkboth\@gobbletwo
        \def\@oddhead{\tt Adler-van\,Moerbeke:%
        Random permutations\hfil Nov 30,1999\hfil\S\thesection, p.\thepage}
        \def\@oddfoot{\rm\hfil\thepage\hfil}
        \let\@evenhead\@oddhead
        \let\@evenfoot\@oddfoot}
\makeatother

\begin{document}
\maketitle
\tableofcontents

%\abstract{ Matrix Fourier-like integrals
%over the classical groups $O_+(n), O_-(n), Sp(n)$ and
%$U(n)$ are connected with the distribution
%of the length of the longest increasing sequence in
%random permutations and random involutions and the
%spectrum of random matrices. One of the purposes of
%this paper is to show that all those integrals satisfy
%the Painlev\'e V equation. In this work, we present
%both, new results and known ones, in a unified way.

%Our method consists of inserting one set of
%time variables
% t=(t_1,t_2,...) in the integrals for the real compact
% groups and two sets of times (t,s) for the
% unitary group. The point is that these new time-dependent
% integrals satisfy integrable hierarchies:
% (i) O(n) and Sp(n) correspond to the
%  standard Toda lattice.
% (ii) U(n) corresponds to the Toeplitz lattice,
%  a very special reduction of the
%  discrete sinh-Gordon equation.

% Both systems, the standard Toda lattice and the Toeplitz
%  lattice are also reductions of the 2-Toda lattice, thus
% leading to a natural vertex operator, and so, a natural
% Virasoro algebra, a subalgebra of which annihilates
% the tau-functions.  Combining these equations leads
 % to the Painlev\'e V equation for the integrals.}

\setcounter{section}{-1} % introduction is section 0
\section{Introduction}

 In recent times, there has been a
 considerable interest in
matrix Fourier-like integrals over the classical
groups $O_+(\ell), O_-(\ell), Sp(\ell)$ and $U(\ell)$,
due to their connection with the distribution of the
length of the longest increasing sequence in random
permutations and random involutions and also with the
spectrum of random matrices. This connection first
appeared in I. Gessel's work \cite{G}, who showed that
some generating function for the distribution of the
length of the longest increasing sequence can be
represented as a Toeplitz matrix. One of the purposes
of this paper is to show that all those expressions
are unique solutions to the Painlev\'e V equation,
with certain initial condition. In this work, we
present both, new results, concerning $O(\ell)$ and
known ones, concerning $U(\ell)$
; all cases are done
in a same unified way.

Our method consists of appropriately adding one set of time
variables
 $t=(t_1,t_2,...)$ to the integrals for the real compact
 groups and two sets of times $(t,s)=(t_1,t_2,...,s_1,s_2,...)$ for the
 unitary group. The point is that these new time-dependent
 integrals satisfy integrable hierarchies:
\begin{description}
  \item[(i)] $O_{\pm}(\ell)$ and $Sp(\ell)$ correspond to the
  {\sl  standard
  Toda lattice}; the associated moment matrices are H\"ankel,
  whose determinants provide the Toda
  $\tau$-functions.
 %their $\tau$-functions are determinants of H\"ankel matrices.

  \item[(ii)] $U(\ell)$ corresponds to
  a very special case of the
  {\sl discrete sinh-Gordon equation}, leading to a
  new lattice, the {\sl Toeplitz lattice}. This
  lattice involves a dual pair of infinite variables
  $x_i$ and $y_i$, themselves matrix integrals.
  Its $\tau$-functions are determinants
  of moment matrices, which are Toeplitz.
\end{description}

 Both systems, the standard Toda lattice and the Toeplitz
  lattice are peculiar reductions of the 2-Toda lattice.
 Each reduction
 has a natural vertex operator, and so, a natural
 Virasoro algebra, a subalgebra of which annihilates
 the $\tau$-functions.  Combining these equations and,
 in the end, evaluating the result along appropriate
 $(t,s)$-loci all lead, in a unifying and quick way, to different versions of
 the Painlev\'e V equation for the integrals.
 More details
 about the precise nature of the Painlev\'e equations will be given in
propositions 3.3,
 4.1 and 4.2. After this paper had been written, we
 found out the Toeplitz lattice coincides with the
 so-called Ablowitz-Ladik system; see Suris \cite{Suris}. However, our
 approach to that system is novel.

Let $S_n$ be the group of $n!$ permutations $\pi_n$
and $S_{2n}^0$ the subset of $(2n-1)!!=
\frac{(2n)!}{2^nn!}$ fixed-point free involutions
$\pi^0$(i.e., $(\pi^0) ^2=I$ and $\pi^0 (k) \neq k$
for $1\leq k \leq 2n$ ). $\pi_n$ refers to
 a permutation in $S_n$ and $\pi^0_{2n}$ to
 an involution in $S^0_{2n}$. Also consider
$S_{n,k}=\{\mbox{words of length $n$ from an alphabet
of $k$ letters}\}$.

An {\em increasing subsequence} of $\pi \in S_n$ or
$S_{n}^0$ is a sequence $1\leq j_1<...< j_k\leq n$,
such that $\pi (j_1)<...<\pi(j_k)$. Define $$
\sigma(\pi_n) =  \mbox{ length of the longest
increasing subsequence of $\pi_n$ }.
 $$
  In the case of $S_{n,k}$, the definition of
$\sigma$ is the same, except that the subsequences
must be increasing, without being necessarily strictly
increasing.

%\noindent
{\sl Notation}: The expectations
$E_{O(\ell)},~E_{U(\ell)},...$ refer to integration
with regard to Haar measure, normalized so that
$E_{O(\ell)}(1)=1,~E_{U(\ell)}(1)=1,...$, as it
should. Sometimes, it will be more convenient to use
integrals $\int_{O(\ell)},~\int_{U(\ell)},...$, which
refer to integration with respect to Haar measure,
normalized as in Proposition 1.1 below. For $U(\ell)$,
the two normalizations happen to agree.

\begin{theorem} For every $\ell\geq 0$,
%\footnote{$dM$ denotes
%Haar measure on the corresponding group normalized
% as in (1.0.2), so that vol$(U(n))=1$. In the
%symplectic
% case below,
%take $\ell$ odd.}
 the generating functions below have the following expression in
terms of specific solutions of the Painlev\'e V
equation
:

  \bea
\mbox{\bf (i)}~~&&
\lefteqn{2\sum^{\iy}_{n=0}\frac{x^{2n}}{(2n)!} \#
\{\pi^0_{2n} \in S^0_{2n} ~\bigl|~ \sigma(\pi^0_{2n})
 \leq \ell +1  \}} \nonumber\\
 &&~~~~~~~~~~~~= E_{O(\ell+1)_{-}}e^{x~tr M}dM
 +  E_{O(\ell+1)_{+}}e^{x~tr M}dM
 %~~\mbox{or}~~
 %\int_{Sp(\frac{\ell-1}{2})}e^{x~tr M}dM
 \nonumber\\
 &&~~~~~~~~~~~~=\exp \left({\int_0^x\frac{f^-_{\ell
 }(u)}{u}du}\right)+\exp \left({\int_0^x\frac{f^+_{\ell
}(u)}{u}du}\right)
%~,~~\mbox{where}~E=\int_{O(\ell+1)_{\pm}}e^{x~tr M}dM~~\mbox{or}~~
%\int_{Sp(\frac{\ell-1}{2})}e^{x~tr M}dM
 \nonumber\\
\mbox{\bf
(ii)}~~&&\sum^{\iy}_{n=0}\frac{x^n}{(n!)^2}\# \{\pi_n
\in S_n ~ \bigl|~ \sigma(\pi_n)\leq\ell \}= E_{U(\ell
)}~e^{\sqrt{x}~tr (M+\bar M)}dM \nonumber\\ &&
 ~~~~~~~~~~~~~~~~~~~~~~~~~~~~~~~~~~~~~~~~=\exp
{\int_0^x \log \left(\frac{x}{u}
\right)g_{\ell}(u)du},
%~,~~\mbox{where}~~
% E=\int_{U(\ell )}e^{\sqrt{x}~tr (M+\bar M)}dM
 \nonumber \\
 \mbox{\bf (iii)}~~&&
 %e^{-\ell x}
 \sum_{n=0}^{\iy}
 \frac{x^{n}}{n!}
\#\left\{ \pi_n \in S_{n,k} ~ \bigl|~
\sigma(\pi_n)\leq \ell \right\}=
 E_{U(\ell)}~ \det (I+M)^k e^{-x~ tr \bar M}dM\nonumber\\
 &&~~~~~~~~~~~~~~~~~~~~~~~~~~~~~~~~~
  ~~~~~~~~
   =\exp \left({x\ell+(\ell+k)\int_0^x\frac{h_{\ell}(u)}{u}du}\right)
 \nonumber \\
\eea
  where $f_{\ell},~ g_{\ell}$ and $h_{\ell}$
are unique solutions to three different versions of
the Painlev\'e V equation, with the initial condition
indicated below; to be precise

 \be
\mbox{\bf (i)} \left\{\begin{array}{l}
\displaystyle{f^{\prime\prime\prime}+
 \frac{1}{u}
f^{\prime\prime}
 +\frac{6}{u}{f^{\prime} }^2-\frac{4}{u^2} f f^{\prime}
 -\frac{16 u^2+\ell^2}{u^2}f^{\prime}
+\frac{16}{u}f  +\frac{2(\ell^2-1)}{u}=0} \\  \\
\displaystyle{\mbox{with}~~  f_{\ell}^{\pm}(u)=u^2 \pm
\frac{u^{\ell+1}}{\ell !}+O(u^{\ell+2}),~\mbox{near}~
u=0 .}
\end{array}
 \right. \ee

\be
\mbox{\bf (ii)} \left\{\begin{array}{l} \displaystyle{
g^{\prime\prime}-
 \frac{g^{\prime 2}}{2}\left(\frac{1}{g-1}+
  \frac{1}{g}\right)+ \frac{g^{\prime}}{u}
   +\frac{2}{u}g(g-1)-
   \frac{\ell^2}{2u^{2}}\frac{g-1}{g}=0}
    \\  \\
\displaystyle{\mbox{with}~~
g_{\ell}(u)=1-\frac{u^{\ell}}{(\ell
!)^2}+O(u^{\ell+1}) ,~\mbox{near}~ u=0 .}
\end{array}
 \right. \ee

\be
\mbox{\bf (iii)} \left\{\begin{array}{l}
\displaystyle{
h^{\prime\prime\prime}-\frac{h^{\prime\prime
2}}{2}\left(
\frac{1}{h'+1}+\frac{1}{h'}\right)+\frac{h''}{u}+\frac{2(\ell
+k)}{u}h'(h'+1)  }\nonumber\\
 \displaystyle{-\frac{1}{{2u^2h^{\prime}(h^{\prime}+1)}}
\Bigl((u-\ell)h^{\prime}-h-\ell \Bigr)
 \Bigl((2h +u+\ell)h^{\prime}+h+\ell \Bigr) =0}\nonumber\\
    \\
\displaystyle{\mbox{with}~~
h_{\ell}(u)=u\frac{k-\ell}{k+\ell}-\frac{u^{\ell +1
}}{(\ell+1)!}\left( k+\ell-1 \atop \ell \right)
  +O(u^{\ell+2}) ,~\mbox{near}~ u=0 .}
\end{array}
 \right. \ee

\end{theorem}

That the orthogonal matrix integrals (i) satisfy
Painlev\'e V is new. The identity (i) involving
orthogonal matrix integrals and random involutions is
due to Rains \cite{R}. That the $U(\ell)$-integral
(ii) satisfies Painlev\'e was first established by
Hisakado \cite{H}, using our methods (see \cite{AvM1})
and then reestablished by Tracy and Widom \cite{TW1},
using methods of functional analysis. The identity
between random permutations and unitary matrix
integrals, via Toeplitz determinants, goes back to
Gessel \cite{G}. Similarly, the $U(\ell)$-integral
(iii) was first established by Tracy-Widom \cite{TW2},
again using methods of functional analysis. The
relation of the combinatorics to integrals over the
groups was extensively studied by Diaconis and
Shashahani \cite{DS}, Rains \cite{R}, Baik and Rains
\cite{BR}; see also Johansson \cite{J}, Baik, Deift
and Johansson \cite{BDJ}, Aldous and Diaconis
\cite{AD}, Tracy and Widom \cite{TW1,TW2}.

 Our methods have the
benefit of providing a unifying (and also quick) way
of establishing these results, new and known ones. The
relationship with integrable systems can be summarized
by Theorems 0.2 and 0.3:

\begin{theorem}

Defining the integrals
   \bean
  {\bf (i)}&& I^{\pm}_{\ell}(x)=\int_{O_{\pm}(\ell)}
e^{x ~tr M }dM  \\
 {\bf (ii)}&&
I_{\ell}(x,y)=\int_{U(\ell)} e^{ tr (xM -y\bar M) }dM,
  \eean
the expressions\footnote{In this statement,
 we use the following notation:
 $$ [n]_{\mbox{\tiny{even}}}:=\max~
\{\mbox{even $x$, such that $x\leq n  \}$.}$$}
 \bean
 {\bf (i)}&& q_{\ell}(x)=\log e_{\ell}^{\pm}~
\frac{I^{\pm}_{\ell+2}}{I^{\pm}_{\ell}},~~\mbox{ with}~~ e^{+}_{\ell}=
 \frac{2}{[\ell+2]_{\mbox{\tiny{even}}}}~~
 \mbox{ and }~e^{-}_{\ell}=
 \frac{2}{[\ell+1]_{\mbox{\tiny{even}}}}\\
 %  ~e^{\pm}_{\ell}=
 %\frac{2}{[\ell+{2 \atop 1}]_{\mbox{\tiny{even}}}}
  {\bf (ii)}&& q_{\ell}(x,y)=\log ~
\frac{I_{\ell+1}}{I_{\ell}}
 \eean
 satisfy respectively
  \bean
{\bf (i)}&& \frac{1}{4}\frac{\pl^2 q_{\ell} }{\pl x^2}  =-e^{q_{\ell}-q_{\ell
-1}}+ e^{q_{\ell +1}-q_{\ell }}~~~~~~~~ \mbox{\bf (standard Toda
lattice)}\\
  {\bf (ii)}&& \frac{\pl^2 q_{\ell}}{\pl
x\pl y}  =e^{q_{\ell}-q_{\ell -1}}- e^{q_{\ell
+1}-q_{\ell }}. ~~~ \mbox{\bf (discrete sinh-Gordon equation)}
 \eean
\end{theorem}

\remark Note, if the lattice is 2-periodic, i.e.,
$q_{\ell }=q_{\ell +2k}$, then (ii) becomes the
sinh-Gordon equation for $r=q_{\ell }-q_{\ell -1}$: $$
\frac{\pl^2r}{\pl x\pl y}=4\mbox{\,sinh\, } r. $$

\vspace*{1cm}

Define the following probability measure on the
unitary group $U(n)$: $$ P^{t,s}_{U(n)}(M\in
dM):=\tau_n(t,s)^{-1}e^{\sum_1^{\iy} \Tr
(t_iM^i-s_i\bar M^i)}dM ,$$
 and
  $ h=\mbox{diag}(h_0,h_1,...),~~h_n=\tau_{n+1}/\tau_n $,
with
 $$ \tau_n(t,s):=\int_{U(n)}e^{\sum_1^{\iy}\Tr
(t_iM^i-s_i\bar M^i)}dM .
 $$

 Also,
let $p_i^{(1)}(t,s;z)$ and $p_i^{(2)}(t,s;z)$ be bi-orthogonal
monic polynomials in $z$, depending on $t$ and $s$, satisfying
$\la p_i^{(1)}(t,s;z),p_j^{(2)}(t,s;z)\ra_{t,s}=\delta_{ij}h_i
$, with regard to the inner-product
$$
\la f(z),g(z)\ra_{t,s}:=\oint_{S^1} \frac{dz}{2\pi i
z}f(z)g(z^{-1})
 e^{\sum_1^{\iy}(t_iz^i-s_iz^{-i})},~~~~t,s\in \BC^{\iy}.
 $$

The statement of Theorem 0.3 contains the elementary
Schur polynomial\footnote{They should not be confused
with the
 bi-orthogonal polynomials $p_i^{(k)}(t,s;z)$. } $p_{n}$,
 defined by $e^{\sum^{\iy}_{1}t_iz^i}:=\sum_{i\geq 0}
 p_i(t_1,t_2,...)z^i$ and applied to the spectrum $x_k=e^{i\theta_k}$
of the unitary matrix $M\in U(n)$:

% $s_{\lambda}(M)$ associated with a Young diagram
%$\lambda$,  In the formula below, the
% $p_{n}$'s denote the elementary Schur polynomials
%$e^{\sum^{\iy}_{1}t_iz^i}:=\sum_{i\geq 0}
% p_i(t_1,t_2,...)z^i$ :
%  $$ (-1)^n
%s_{\lambda}(M):=p_n(-\sum_k x_k,
% -\frac{1}{2}\sum_k x^2_k,-\frac{1}{3}\sum_k x^3_k,...)
%~~\mbox{for}~~\lambda=({\underbrace{1,...,1}_n}).
% $$

\begin{theorem} %Given the Young diagram
%$ \lambda=(~\overbrace{1,...,1}^n~)$,
Consider the following
 variables, expressed in terms of the expectation
 for the distribution above, or expressed in terms
 of the bi-orthogonal polynomials evaluated at $z=0$,
 \bea
 x_n (t,s):&=& E^{t,s}_{U(n)}p_n(-\Tr
 M,-\frac{1}{2}\Tr M^2,-\frac{1}{3}\Tr M^3,...)
% \frac{ \int_{U(n)} p_n(-Tr
% M,-\frac{1}{2}Tr M^2,-\frac{1}{3}Tr M^3,...)
%s_{\lambda} (M)
%e^{\sum_1^{\iy} Tr (t_iM^i-s_i\bar
%M^i)} dM}{(-1)^n \int_{U(n)}  e^{\sum_1^{\iy} Tr
%(t_iM^i-s_i\bar M^i)}dM}
\nonumber\\
&=&\frac{p_n(-\tilde
\pl_t)\tau_n(t,s)}{\tau_n(t,s)}=p_n^{(1)}(t,s;0)
\nonumber\\
 y_n(t,s):&=& E^{t,s}_{U(n)}p_n (-\Tr \bar M,
 -\frac{1}{2}\Tr \bar M^2,-\frac{1}{3}\Tr \bar M^3,...)
% \frac{ \int_{U(n)}
%  p_n(-Tr \bar M,-\frac{1}{2}Tr \bar M^2,
%  -\frac{1}{3}Tr \bar M^3,...)
%s_{\lambda} (\bar M)
%e^{\sum_1^{\iy} Tr
%(t_iM^i-s_i\bar M^i)} dM}{(-1)^n \int_{U(n)}
%e^{\sum_1^{\iy} Tr (t_iM^i-s_i\bar M^i)}dM}
\nonumber
\\ &=&\frac{p_n(\tilde \pl_s)\tau_n(t,s)}{\tau_n(t,s)}
 =p_n^{(2)}(t,s;0).
  \eea
%& & $~~~~~~~~~\displaystyle{=\frac{p_n(\tilde
%\pl_s)\tau_n(t,s)}{\tau_n(t,s)}}$\\
%\end{tabular}}
The $x_n$ and $y_n$'s satisfy the following integrable
Hamiltonian system
 \bea
 \frac{\pl x_n}{\pl
t_i}=(1-x_ny_n)\frac{\pl H^{(1)}_i}{\pl y_n}  &~~~~~&
\frac{\pl y_n}{\pl t_i}=-(1-x_ny_n)\frac{\pl
H^{(1)}_i}{\pl x_n} \nonumber \\ \frac{\pl x_n}{\pl
s_i}=(1-x_ny_n)\frac{\pl H^{(2)}_i}{\pl y_n}  &~~~~~&
\frac{\pl y_n}{\pl s_i}=-(1-x_ny_n)\frac{\pl
H^{(2)}_i}{\pl x_n}
  , \\ &&\hspace{2.5cm}(\mbox{\bf Toeplitz lattice})
 \nonumber \eea
%for the symplectic structure $$ \omega := \sum_1^{\iy}
%\frac{dx_k \wedge dy_k}{1-x_ky_k}, $$ and the
with initial condition $x_n(0,0)=y_n(0,0)=0$ for
$n\geq 1$ and boundary condition
$x_0(t,s)=y_0(t,s)=1$. The traces
 $$ H^{(k)}_i=-
\frac{1}{i}\Tr~{ L} _k^i,~~i=1,2,3,...,~~k=1,2. $$ of
the matrices ${ L}_i$ below are integrals in
involution with regard to the symplectic structure
 $$ \omega := \sum_1^{\iy}
\frac{dx_k \wedge dy_k}{1-x_ky_k}, $$
  where $L_1$ and $L_2$ are given by the `` rank 2" semi-infinite
matrices $$h^{-1}L_1 h:= \left(\begin{tabular}{lllll}
$-x_1y_0$  &  $1-x_1y_1$ & ~~ $0$      & ~~ $0$ &   \\
$-x_2y_0$ &  $-x_2y_1$ & $1-x_2y_2$& ~~ $0$   & \\
$-x_3y_0$ &  $-x_3y_1$ & $ -x_3y_2$&  $1-x_3y_3$ & \\
$ -x_4y_0$ &  $ -x_4y_1$ & $-x_4y_2$  & $ -x_4y_3$   &
\\
 & &  &    &  $\ddots$\\
\end{tabular}
\right)
$$
and
$$L_2:=
\left(\begin{tabular}{lllll}
$-x_0y_1$  &  $-x_0y_2$ & $-x_0y_3$     & $-x_0y_4$ &   \\
$1 -x_1y_1$ &  $-x_1y_2$  & $-x_1y_3$& $-x_1y_4$   & \\
~~$0$       &  $1 -x_2y_2$ & $ -x_2y_3$&  $-x_2y_4$ & \\
~~$0$       &  ~~$0$      & $ 1 -x_3y_3$  & $ -x_3y_4$   &  \\
 & &  &    &  $\ddots$\\
\end{tabular}
\right). $$ Morover, the precise ``rank 2"-structure
of $L_1$ and $L_2$ is preserved by the equations
 $$ \frac{\pl
L_i}{\pl t_n}=\bigl[\bigl(L_1^n\bigr)_+,L_i\bigr]
\quad\hbox{and}\quad \frac{\pl L_i}{\pl
s_n}=\bigl[\bigl(L_2^n\bigr)_-,L_i\bigr] \quad
i=1,2~\mbox{and}~n=1,2,\dots
$$
\hfill ({\bf Two-Toda Lattice})

\end{theorem}

\remark The first equation in the hierarchy above,
corresponding to the Hamiltonians
 $$
 H_1^{(1)}=-\Tr~ L_1=\sum_0^{\iy} x_{i+1}y_i,~
 H_1^{(2)}=-\Tr~ L_2=\sum_0^{\iy} x_{i}y_{i+1},
 $$
 reads: \bean
\frac{\pl x_n}{\pl t_1}=x_{n+1}(1-x_ny_n) &~~~~~&
\frac{\pl y_n}{\pl t_1}=-y_{n-1}(1-x_ny_n)  \\
\frac{\pl x_n}{\pl s_1}=x_{n-1}(1-x_ny_n)  &~~~~~&
\frac{\pl y_n}{\pl s_1}=-y_{n+1}(1-x_ny_n) . \\ \eean

\newpage

Here, we {\bf outline the ideas and the results} in
the paper. Throughout, consider a weight $\rho(x)dx$
on an interval $F\subset\BR$,
satisfying\footnote{Decaying rapidly means:
$\rho(x)f(x)=0$ at finite boundary points of $F$, or
$\rho(x) f(x)x^k\rightarrow 0$, when $x\rightarrow
\{\mbox{an infinite boundary point}\}$, for all
$k=0,1,2,\ldots$\,. }
\be
-\frac{\rho'(x)}{\rho(x)}=\frac{\sum_{i\geq 0} b_i x^i}
{\sum_{i\geq 0} a_i x^i}=\frac{g(x)}{f(x)},~~\mbox{with $\rho(x)$
decaying rapidly at $\pl F$}.
\ee
We now define two time-dependent inner-products, one
 given by a weight $\rho(x)dx$ on the real line $\BR$
 and another given by a contour integration about the unit circle
 $S^1\subset \BC$,
\be
\left\{\begin{array}{l}
\displaystyle{\la f(x),g(x)\ra_t:=\int_{\BR} f(x)g(x)
 e^{\sum_1^{\iy}t_ix^i}
 \rho(x) dx},~~t\in \BC^{\iy}
 %\nonumber
 \\
 \displaystyle{ \la f(z),g(z)\ra_{t,s}:=\oint_{S^1} \frac{dz}{2\pi i
z}f(z)g(z^{-1})
 e^{\sum_1^{\iy}(t_iz^i-s_iz^{-i})},~~~~t,s\in \BC^{\iy}
 }\end{array}
\right. .
\ee
These inner-products lead to H\"ankel and Toeplitz moment matrices,
respectively,
 \be
\left\{\begin{array}{l}
 \displaystyle{m_n(t):=\left(\la x^i,x^j\ra_{t}
 \right)_{0\leq i,j\leq n-1} }~~~~~~
 \mbox{({\bf H\"ankel})}
 \nonumber\\  \nonumber\\
 \displaystyle{m_n(t,s):=\left(\la z^i,z^j\ra_{t,s}
 \right)_{0\leq i,j\leq n-1}}
 ~~\mbox{({\bf Toeplitz}).}
 \end{array}
\right.
\ee
  The determinants $\tau_n$ of the $m_n$'s
have different representations: on the
 one hand, as multiple integrals, involving Vandermonde's
 $\Delta_n(z)$, and, on the other hand,
 as inductive expressions in term of $\tau_{n-1}$,
 involving a vertex operator\footnote{For
 $v=(v_{0},v_{1},\ldots)^{\top},~
  (\Lambda v)_{n}=v_{n+1},
  ~(\Lambda^{\top}v)_{n}=v_{n-1}$, and $\chi(z):=(1,z,z^2,...)$.}
 \be
\BX_{12}(t,s;u,v)=\Lb^{\top}e^{\sum_1^{\iy}(t_i u^i-s_i v^{i})}
e^{-\sum_1^{\iy}(\frac{ u^{-i}}{i} \frac{\pl}{\pl t_i}-
 \frac{ v^{-i}}{i}
\frac{\pl}{\pl s_i})}\chi(uv),
 \ee
  to be explained in (0.0.9). The
$\tau_n(t)$ and $\tau_n(t,s)$
 are respectively solutions to the standard Toda lattice,
and the so-called Toeplitz lattice, both reductions of
the semi-infinite 2d-Toda lattice\footnote{The
expression
  $\BX_{12}(\frac{s+t}{2},\frac{s-t}{2};u,u)$ is actually independent
  of $s$ ! The expressions in (0.0.9) are inductive in $n$,
  because of the presence of the downwards shift $\Lb^{\top}$ in $\BX_{12}$.   }:
 \bea
\hspace{-0.1cm}I_n&=&n!\det \tau_n=n!\det m_n \nonumber\\
\hspace{-0.3cm}&=&\left\{\begin{array}{l}
\displaystyle{\int_{\BR^n} \Delta_n^2\prod_{k=1}^n
  e^{\sum_{i=1}^{\iy} t_i z_k^i} \rho(z_k) dz_k}
 =\int_{\BR}du \rho(u) \left(\BX_{12}\bigl(\frac{s+t}{2},
  \frac{s-t}{2};u,u\bigr)
 I \right)_{n} \nonumber\\
 \hspace{6cm}(\mbox{\bf standard Toda $\tau$-functions})\nonumber\\
\displaystyle{
 \oint_{(S^1)^{n}}|\Dt_n|^{2}
 \prod_{k=1}^n
e^{\sum_1^{\iy}(t_i z_k^i-s_iz_k^{-i})}
 \frac{dz_k}{2\pi i z_k}=\int_{S^1} \frac{du}{2\pi i u}
\left(\BX_{12}(t,s;u,u^{-1})I \right)_{n},
 }%\nonumber
 \\
 \hspace{6cm}(\mbox{\bf two-Toda $\tau$-functions})\nonumber\\
 \end{array}
\right.\\
\eea
where $\tau_n(t)$ and $\tau_n(t,s)$ satisfy  the following differential
equations, (the second one is new)
\bea
\left\{
\begin{array}{l}
\displaystyle{\frac{\pl^4}{\pl t_1^4}\log\tau_n
+6\left(\frac{\pl^2}{\pl t_1^2}\log\tau_n\right)^2
+3\frac{\pl^2}{\pl t_2^2}\log\tau_n
-4\frac{\pl^2}{\pl t_1\pl t_3}\log\tau_n = 0,}\nonumber\\
\hspace{6cm}(\mbox{\bf KP}-equation)\nonumber\\
\displaystyle{\frac{\pl^2}{\pl s_2\pl
t_1}\log\tau_{n}=-2\frac{\pl}{\pl
s_1}\log\frac{\tau_{n}}{\tau_{n-1}}~.\frac{\pl^{2}}{\pl s_{1}\pl
t_{1}}\log\tau_{n}- \frac{\pl^{3}}{\pl s_{1}^2\pl
t_{1}}\log\tau_{n}.}\\
\hspace{7cm}\mbox{({\bf
two-Toda}-equation)}\nonumber
\end{array}
\right. \eea The unique factorization of the
time-dependent semi-infinite moment matrices
$m_{\infty}$ (defined just under (0.0.7)) into lower-
times upper-triangular matrices
\be
\left\{
\begin{array}{l}
m_{\infty}(t)=S(t)^{-1}S^{\top -1}(t)\\
m_{\infty}(t,s)=S_{1}(t,s)^{-1}S_2(t,s)
\end{array}
\right. \ee leads to matrices $S(t)$,  $S_{1}(t,s)$,
$S_{2}(t,s)$  of the form $$S(t)=\sum_{ i\leq
0}a_i\Lambda^i,\quad S_{1}(t,s)=\sum_{i\leq
0}b_i\Lambda^i,\quad S_{2}(t,s)=\sum _{i\geq 0}b'_i
\Lambda^i $$ with $b_{0}=I$, $a_{i},b_{i},b_{i}'$
diagonal matrices. By ``dressing up" the shift $\Lb$ (
defined in footnote 4), they evolve according to the
following integrable systems:

  \medbreak
$$\left\{
 \begin{array}{l}  L(t)=S(t)\Lambda
S^{-1}(t) :\mbox{symmetric and tridiagonal, ~~~ ({\bf
Toda lattice})}
 \\  \\
 \displaystyle{\left\{
\begin{array}{l} L_1(t,s)=S_{1}(t,s)
\Lambda S_{1}^{-1}(t,s)
  \\
L_2(t,s)=S_{2}(t,s)\Lambda^{\top} S_{2}^{-1}(t,s).
\end{array} ~~~~~~~~~~~~~~~~~~~~~~~~~
~\mbox{({\bf Toeplitz lattice})}
\right. }
\end{array} \right. $$

%\vspace{1cm}

\noindent As already pointed out in (0.0.9), that
expression involves the following reduction of the
two-Toda vertex operator $\BX_{12}$:
 \bea \left\{\begin{array}{l}
 \displaystyle{\BX_{12}\bigl(\frac{s+t}{2},\frac{s-t}{2};u,u\bigr)
 =:\BX(t;u)=\Lambda^{\top}\chi(u^2)
 e^{\sum_1^{\infty} t_iu^i}e^{-2\sum_1^{\infty}
\frac{u^{-i}}{i} \frac{\partial}{\partial t_i}}
}
 \\  \\
 \displaystyle{\BX_{12}(t,s;u,u^{-1})=\Lb^{\top}
  e^{\sum_1^{\iy}(t_i u^i-s_i u^{-i})}
e^{-\sum_1^{\iy}(\frac{ u^{-i}}{i}
\frac{\pl}{\pl t_i}-
 \frac{ u^{i}}{i}
\frac{\pl}{\pl s_i})}}
 \end{array}
\right.
\eea

Each of these vertex operators leads to Virasoro
algebras $\JR^{(2)}_{m}$ and $\VR^{(2)}_{m}$ of
central charge $c=1$ and $c=0$ respectively, defined
by  \bea \left\{\begin{array}{l}
 \displaystyle{\frac{\pl}{\pl u}u^{m+1}f(u)
  \BX_{}(t,u)\rho(u)
 =
\Bigl[ \JR^{(2)}_m(t),
\BX_{}(t,u)\rho(u)\Bigr]}\nonumber\\
 \nonumber\\
 \displaystyle{\frac{\pl}{\pl u}u^{m+1} \frac{\BX_{12}(t,s;u,u^{-1})}{u}
 =\left[ {\cal V}^{(2)}_m(t,s),
  \frac{\BX_{12}(t,s;u,u^{-1})}{u}\right],}
 \end{array}
\right. \eea and having the explicit expressions (see
notation (0.0.6))

\bea
\JR^{(2)}_m(t)&:=&\sum_{i\geq 0}\left.\left( a_i \,^{\beta}\BJ_{i+m}
^{(2)}(t)-
b_i \,^{\beta}\BJ_{i+m+1}^{(1)}(t)\right)\right|_{\beta =2}\nonumber\\
 {\cal V}^{(2)}_m(t,s)&:=&\left.{}^{\beta}\BJ_{m}^{(2)}(t)-
  {}^{\beta} \BJ_{-m}^{(2)}(-s)
-m\left( \theta~~ {}^{\beta}\BJ_{m}^{(1)}(t)+(1-\theta)~~
 {}^{\beta}\BJ_{-m}^{(1)}(-s) \right)\right|_{\beta=1}\nonumber\\
 & &
 \eea
in terms of generators $^{\beta}\BJ_{m}^{(2)}$ defined in (5.0.3)
below, and arbitrary $\theta$.

The point is that a big subalgebra of $\JR_{m}^{(2)}$'s  and a
small one of $\VR_{m}^{(2)}$'s annihilate $\tau_{n}(t)$ and
$\tau_{n}(t,s)$ respectively, for appropriate $\theta$ and for all
$n\geq 0$,
\be
\left\{\begin{array}{l}
 \displaystyle{\JR^{(2)}_m\tau_n(t)=0~~\mbox{for}~~ m\geq -1,}
 \nonumber\\  \nonumber\\
 \displaystyle{{\cal V}^{(2)}_m\tau_n(t,s)=0~~\mbox{for}~~ m=-1,0,1~~
  \mbox{($SL(2,\BZ)$- algebra).}
}
 \end{array}
\right.
\ee

\noindent To summarize, we have that combining
 these equations and restricting to the three
 different loci $\LR$ below, always leads to Painlev\'e V:

$$
\left.\left\{
\begin{array}{l}
\mbox{KP\,} (\tau_n)=0\\
\JR^{(2)}_m\tau_n(t)=0,\\ ~~\mbox{for}~~m=-1,0
\end{array}
\right\}\right|_{\LR=\left\{
\begin{array}{l}
t_{1}=x, \mbox{\,all other}\\
t_{i}=0
\end{array}
\right\}}\Longrightarrow\left\{\begin{array}{l} \mbox{Painlev\'e $V$
for}\\ O_{\pm}(n)\mbox{-integral}%~\mbox{or} ~Sp(n)  \\
\\
\end{array}
\right.
$$

\vspace*{1cm}
$$
\left.\left\{
\begin{array}{l}
2-\mbox{Toda PDE}\\
\VR^{(2)}_{m}\tau_n(t,s)=0,\\~~\mbox{for}~~ m=-1,0,1 \\
\mbox{Toeplitz relation}
\end{array}
\right\}\right|_{{\LR=\left\{
\begin{array}{l}
t_{1},s_1\neq 0, \mbox{\,all other}\\
\mbox{all other\,} t_{i},s_i=0
\end{array}
\right\}~\mbox{or}} \atop{\LR=\left\{
\begin{array}{l}
\mbox{all}~ it_i=-k(-1)^i \\
\mbox{$s_i=0$, except\,} s_{1}=x
\end{array}
\right\}}}
\Rightarrow
\left\{\begin{array}{l}
\mbox{Painlev\'e $V$}\\ \mbox{for} \\
U(n)\mbox{-integral}
\end{array}  \right.
$$

Acknowledgement: We wish to thank Jinho Baik, Craig
Tracy and Harold Widom for several informative
discussions.

\section{Integrals over classical groups and combinatorics  }

%\subsection{Integrals over classical groups and Tchebychev polynomials}

%$^{\beta} G$

This section contains a number of useful facts about
integrals over groups, its relation with combinatorics
and finally the behavior of some of the integrals near
$x=0$.

 The situation is quite different, according to
whether one integrates over the real $(O_{\pm}$,
$Sp(\ell))$ or the complex $(U(\ell))$. The real group
integrals involve the
 Jacobi weight,
\be
\rho_{\alpha\beta}(z)dz:=(1-z)^{\alpha}(1+z)^{\beta} dz,
\ee
for $\alpha,~\beta=\pm 1/2$ and the Tchebychev
 polynomials $T_n(z)$, defined by $T_n(cos \theta)\linebreak :=cos n \theta$.
 In particular, we have $T_1(z)=z$. We now have the following theorem (see
Johansson \cite{J}):

\begin{proposition} ({\bf Weyl})
 Defining $$g(z):=  2\sum_1^{\iy}t_iT_i(z),$$ the following holds:
\bea
\int_{U(n)}e^{\sum_1^{\iy}tr (t_iM^i-s_i\bar
M^i)} dM &=&
 \frac{1}{n!}
 \int_{(S^1)^{n}}|\Dt_n(z)|^{2}
 \prod_{k=1}^n
e^{\sum_1^{\iy}(t_i z_k^i-s_iz_k^{-i})}
 \frac{dz_k}{2\pi i z_k}
\nonumber\\ \int_{O(2n+1)_+}e^{\sum_{1}^{\iy} t_i \Tr
M^i}dM&=& e^{\sum_{1}^{\iy} t_i}
 \int_{[-1,1]^n}\Delta_n(z)^2
 \prod_{k=1}^n e^{g(z_k)}
  \rho_{(\frac{1}{2},-\frac{1}{2})}(z_k) dz_k\nonumber\\
  \int_{O(2n+1)_-}e^{\sum_{1}^{\iy} t_i \Tr M^i}dM
  &=& e^{\sum_{1}^{\iy} (-1)^i t_i}
 \int_{[-1,1]^n}\Delta_n(z)^2
 \prod_{k=1}^n e^{g(z_k)}
  \rho_{(-\frac{1}{2},\frac{1}{2})}(z_k) dz_k\nonumber\\
  \int_{O(2n)_+}e^{\sum_{1}^{\iy} t_i t\Tr M^i}dM
  &=&
 \int_{[-1,1]^n}\Delta_n(z)^2
 \prod_{k=1}^n e^{g(z_k)}
  \rho_{(-\frac{1}{2},-\frac{1}{2})}(z_k) dz_k\nonumber\\
  \int_{O(2n)_-}e^{\sum_{1}^{\iy} t_i \Tr M^i}dM
  &=& e^{\sum_{1}^{\iy} 2 t_{2i}}
 \int_{[-1,1]^{n-1}}\Delta_{n-1}(z)^2
 \prod_{k=1}^{n-1} e^{g(z_k)}
  \rho_{(\frac{1}{2},\frac{1}{2})}(z_k) dz_k\nonumber\\
   \int_{Sp(n)}e^{\sum_{1}^{\iy} t_i \Tr M^i}dM
  &=&
 \int_{[-1,1]^n}\Delta_n(z)^2
 \prod_{k=1}^n e^{g(z_k)}
  \rho_{(\frac{1}{2},\frac{1}{2})}(z_k) dz_k.\nonumber\\
 \eea
 With this normalization, we have (see Appendix 3)
  \bean
  \int_{U(n)}dM&=&1 \\
  \int_{O(2n+1)_{\pm}}dM&=&2^{n^2}
 \displaystyle{\prod_{j=1}^n
   \frac{j!(j-1/2)\Gamma^2(j-1/2)}{(n+j-1)!}} \\
  \int_{O(2n)_{+}}dM&=&\displaystyle{2^{n(n-1)}\prod_{j=1}^n
   \frac{j! \Gamma^2(j-1/2)}{(n+j-2)!}}\\
  \int_{O(2n)_{-}}dM&=&\displaystyle{2^{n(n-1)}\prod_{j=1}^{n-1}
   \frac{j! \Gamma^2(j+1/2)}{(n+j-1)!}}.\eean

\end{proposition}

Letting $\iota \in S_n$ denote the permutation $k
\rightarrow n+1-k$, we also state:

\begin{proposition} The combinatorial quantities
 below have an expression in terms of integrals
 over groups:
  \bean
 \sum_{n\geq 0}
 \frac{x^{2n}}{(n!)^2}
 \#\{\pi \in S_{n},~\sigma_{n}(\pi)\leq \ell  \} &=&
  E_{U(\ell)}e^{x \Tr (M+\bar M)}     \\
 \sum_{n\geq 0}
 \frac{x^{2n}}{(2n)!}
 \#\left\{
\begin{array}{l}
\pi \in S_{4n},~\pi^2=1,~(\pi \iota)^2=1
  \\
\pi(y)
  \neq y,\iota y,~\sigma_{4n}(\pi)
  \leq 2\ell
\end{array}
\right\}
   &=&
  E_{U(\ell)}e^{x \Tr (M+\bar M)}     \\
2 \sum_{n\geq 0}
 \frac{x^{2n}}{(2n)!}
\#\left\{
\begin{array}{l}
\pi \in S_{2n},~\pi^2=1,~
  \\
\pi(y)\neq y,~\sigma_{2n}(\pi)
  \leq \ell
\end{array}
\right\}
  &=& E_{O_-(\ell )}e^{x \Tr M}+E_{O_+(\ell )}e^{x \Tr M}
     \\
\sum_{n\geq 0}
 \frac{x^{n}}{n!}
\#\left\{ \pi \in S_{n},~\pi^2=1,~\sigma_n(\pi)\leq
\ell \right\} &=&e^{x}E_{O_-(\ell+1 )}e^{x \Tr M}
 \\
 \sum_{n\geq 0}
 \frac{x^{n}}{n!}
\#\left\{ \pi \in S_{n},~(\iota
\pi)^2=1,~\sigma_n(\pi)\leq \ell \right\}
&=&e^{x}E_{O_-(\ell+1 )}e^{x \Tr M}
 \\
   \sum_{n\geq 0}
 \frac{x^{2n}}{(2n)!}
 \#\left\{\begin{array}{l}
\pi \in S_{2n},~(\pi \iota)^2=1,
  \\
~\pi(y)\neq \iota y,
  ~\sigma_{2n}(\pi)
  \leq 2\ell
\end{array}\right\}
 &=& E_{O_-(2\ell +2)}e^{x \Tr M}    \\
  \sum_{n\geq 0}
 \frac{x^{2n}}{(n!)^2}
 \#\{\pi \in S_{2n},~\pi \iota=\iota \pi,
  ~\sigma_{2n}(\pi)
  \leq 2\ell  \}  &=&
   E^2_{U(\ell)}e^{x \Tr (M+\bar M)}
     \\
\sum_{n\geq 0}
 \frac{x^{2n}}{(n!)^2}
\#\left\{\begin{array}{l} \pi \in S_{2n},~\pi
\iota=\iota \pi
  \\
\sigma_{2n}(\pi)
  \leq 2\ell +1
\end{array}\right\}
&=&
   E_{U(\ell)}e^{x \Tr (M+\bar M)}  \\
&& \hspace{1cm}
   .E_{U(\ell +1)}e^{x \Tr (M+\bar M)}
     \\
\sum_{n\geq 0}
 \frac{x^{n}}{n!}
\#\left\{ \pi \in S_{n,k},~\sigma_n(\pi)\leq \ell
\right\} &=&E_{U(\ell)}\det (I+M)^k e^{x \Tr \bar M}.
\eean
\end{proposition}

The first identity goes back to Gessel \cite{G} and in
this precise form to Rains \cite{R}. The third,
seventh and eight equalities between the combinatorics
and the integral expression above are due to Rains
\cite{R}, the next ones are due to Baik and Rains
\cite{BR} and the last one is due to Tracy and Widom
\cite{TW2}.

 We now state the following elementary
lemma:

\vspace{1cm}

\begin{lemma}
\bean
 e^{x^2/2}&=&\sum_{n=0}^{\iy} \frac{x^{2n}}{(2n)!}\# \{ \pi
 \in S_{2n}~\bigl|~\pi^2=1,~\mbox{fixed-point
 free}\}=\sum_{n=0}^{\iy}  \frac{x^{2n}}{(2n)!}
 \frac{(2n)!}{2^n n!}\\
 e^{x^2/2+x}&=&\sum_{n=0}^{\iy} \frac{x^{n}}{(n)!}\#\{\pi
 \in S_{n}~\bigl|~\pi^2=1\}=\sum_{n=0}^{\iy}
 \frac{x^{n}}{(n)!} \sum_{0\leq m \leq [n/2]}\left( {n
 \atop 2m} \right) \frac{(2m)!}{2^m m!}\\
 \eean

\end{lemma}

We now estimate the following integrals over
$O(\ell+1)$ and $U(\ell)$ near $x=0$. In all cases,
one notices a big gap in the expansion, -roughly
speaking- of the order $\ell$. This would be hard to
obtain at the level of the integrals, but easy to
obtain via combinatorics.

\begin{proposition} The following estimates hold, near
$x=0$,

\bea
% \lefteqn{e^{-x}\sum_{n=0}^{\iy}
%\frac{x^{n}}{(n)!}\#\{\pi
% \in S_{n}~\bigl|~\pi^2=1,~\sigma_n(\pi)\leq \ell\}  }\\
 E_{O_{\pm}(\ell  +1)} e^{x \Tr M}  &=&
\exp \left(
 \frac{x^2}{2}\pm\frac{x^{\ell+1}}{(\ell+1)!}
 +O(x^{\ell+2})\right),\nonumber\\
 E_{O_+(\ell  +1)} e^{x \Tr M}+E_{O_-(\ell  +1)} e^{x \Tr M}  &=&
 2\exp \left(
 \frac{x^2}{2}
 +O(x^{\ell+2})\right).
\eea

\end{proposition}

\proof From the second relation in Lemma 1.3, it
follows that
 \bean
\lefteqn{ \#\{\pi
 \in S_{n}~\bigl|~\pi^2=1,~\sigma_n(\pi)\leq \ell\}
 }\\
 && \hspace{1cm}= \#\{\pi
 \in S_{n}~\bigl|~\pi^2=1 \} = \sum_{0\leq m \leq
 [n/2]}\left ( {n \atop 2m}\right) \frac{(2m)!}{2^m
 m!},  \mbox{ for } n\leq \ell  \\ &&\hspace{1cm}=  \#\{\pi
 \in S_{\ell+1}~\bigl|~\pi^2=1 \}-1,  \mbox{ for } n= \ell+1.\\
 \eean
 Hence we have from the fourth identity of
 proposition 1.2, and Lemma 1.3,
 \bean
  e^x ~ E_{O_-(\ell+1)}  e^{x \Tr
 M}&=&\sum_{n\geq 0}
 \frac{x^{n}}{n!}
\#\left\{ \pi \in S_{n},~\pi^2=1,~\sigma_n(\pi)\leq
\ell \right\} \nonumber\\ &=&\exp \left(
\frac{x^2}{2}+x -\frac{x^{\ell+1}}{(\ell+1)!}
 +O(x^{\ell+2})\right),
 \eean
  and so
 \be E_{O_-(\ell+1)}  e^{x \Tr
 M}dM=\exp \left( \frac{x^2}{2}
-\frac{x^{\ell+1}}{(\ell+1)!}
 +O(x^{\ell+2})\right).\ee
 But, for
 $2n \leq \ell+1$,
\bean \hspace{-2cm}\lefteqn{ \# \{ \pi
 \in S_{2n}~\bigl|~\pi^2=1,~\sigma_{2n}(\pi)\leq \ell +1,
  ~\mbox{fixed-point
 free}\}  }\\
  &&\hspace{.5cm}=\{\pi \in S_{2n},~\pi^2=1,~ \mbox{fixed-point
 free}\}\\
 &&\hspace{.5cm}= \frac{(2n)!}{2^n n!}, \eean
 and so, from the third identity of Proposition 1.2,
 \bean
\lefteqn{ E_{O_-(\ell+1)}  e^{x \Tr
 M}+ E_{O_+(\ell+1)} e^{x \Tr
 M}}\\
  &=& 2\sum_{n=0}^{\iy} \frac{x^{2n}}{(2n)!}\# \{
\pi
 \in S_{2n}~\bigl|~\pi^2=1,~\sigma_{2n}(\pi)\leq \ell +1,
  ~\mbox{fixed-point
 free}\}\\
 &=& 2 \exp \left(x^2/2 +O(x^{\ell+2})\right),
\eean
 establishing the second relation (1.0.3). Combining
  this formula with (1.0.4), leads to the following
 estimate, near $x=0$
 $$E_{O_{+}(\ell+1)}  e^{x \Tr
 M}
 %= \exp \int_0^x\frac{f_{\ell}^{\pm}(u)}{u} du
   = \exp \left(\frac{x^2}{2} + \frac{x^{\ell+1}}{(\ell+1)!}
    +O(x^{\ell+2})\right),
 $$
establishing the first relation (1.0.3). \qed

%$$ \int_{O_{+}(\ell+1)}  e^{x ~tr
% M} dM =
% \exp \int_0^x\frac{f_{\ell}^{\pm}(u)}{u} du
% =
%  \exp \left(\frac{x^2}{2} \pm \frac{x^{\ell+1}}{(\ell+1)!}
%    +O(x^{\ell+2})\right)$$

%$$
% f_{\ell}^{\pm}(x)=x^2\pm \frac{x^{\ell +1}}{\ell
% !}+O(x^{\ell+1})
% $$

\begin{proposition} The following estimates hold, near
$x=0$,
 \bean E_{U(\ell)}e^{\sqrt{x} \Tr (M+\bar M)}
d M &=&\exp\left(x-\frac{x^{\ell+1}}{((\ell+1)!)^{2}}
  +O(x^{\ell+2})\right)  \\
 E_{U(\ell)}\det(I+M)^k e^{-x\Tr \bar
M} d M &=&
 \exp\left(kx- \frac{x^{\ell+1}}{(\ell+1)!}
  \left( k+\ell \atop \ell +1
\right)
  +O(x^{\ell+2})\right).
\eean

\end{proposition}

\proof Using the first identity of proposition 1.2, we
have
  \bean
I_{\ell}&:=& E_{U(\ell)}e^{\sqrt{x} ~\Tr (M+\bar M)} d
M
\\ &=&\sum_0^{\iy} \frac{x^n}{(n!)^2}\#\{ \pi \in
S_n~\bigl|~\sigma_n(\pi)\leq \ell\}\\
  &=& \sum_0^{\ell} \frac{x^n}{n!}+
  \frac{x^{\ell+1}}{(\ell+1)!^2}
   ((\ell+1)!-1)  +O(x^{\ell+2})\\
  &=&\exp\left(x-\frac{x^{\ell+1}}{(\ell+1)!^2}
  +O(x^{\ell+2})\right).
 \eean

%It follows that $$ g_{\ell}(x)=\frac{d}{dx} x
%\frac{d}{dx}\log I_{\ell}=1  - \frac{x^{\ell}}{(\ell
%!)^2}+O(x^{\ell+1}).$$

%@@@@@@@@@@@@@

 Since the number of (increasing) sequences
$ \overbrace{(1,1,...,2,2,2,...,k,k,k)}^{\ell +1}$  of
length $\ell+1$ and consisting of $k$ symbols, is
given by $\left( k+\ell \atop \ell+1 \right)$, one
computes for $S_{\ell,k}=\{\mbox{words of length
$\ell$ from an alphabet of $k$ letters}\}$:

 $$\#\{ \pi \in
S_{n,k}~\bigl|~\sigma_{n}(\pi)\leq \ell\}=k^{n},
~~\mbox{for}~n\leq \ell,$$

$$\#\{ \pi \in
S_{\ell+1,k}~\bigl|~\sigma_{\ell+1}(\pi)\leq \ell\} =
k^{\ell+1}-\left(k+\ell \atop \ell+1  \right),~~
 \mbox{ for } n= \ell+1.\\
 $$
 Therefore, using the last identity of Proposition 1.2,
one finds
  \bean
I_{\ell}&:=& E_{U(\ell)}\det(I+M)^k e^{-x\Tr \bar M} d
M
\\ &=&\sum_0^{\iy} \frac{x^n}{n!}\#\{ \pi \in
S_{n,k}~\bigl|~\sigma_n(\pi)\leq \ell\}\\
  &=& \sum_0^{\ell } \frac{x^n}{n!}k^n+
   \frac{x^{\ell+1}}{(\ell+1)!}\left(k^{\ell+1}-
   \left(k+\ell \atop \ell+1  \right)\right)
     +O(x^{\ell+2})\\
  &=&\exp\left(kx- \frac{x^{\ell+1}}{(\ell+1)!}
   \left(k+\ell \atop \ell+1  \right)
  +O(x^{\ell+2})\right).
 \eean  \qed

%It follows that $$
%h_{\ell}(x)=\frac{1}{\ell+k}(x\frac{d}{dx} \log
%I_{\ell} -\ell x) =x\frac{k-\ell}{k+\ell}  -
%x^{\ell+1}\left( k+\ell-1 \atop \ell
%\right)+O(x^{\ell+2}).$$

\section{Two-Toda lattice and reductions (H\"ankel and Toeplitz)}

\subsection{Two-Toda on Moment Matrices
and Identities for $\tau$-Functions}

Two-Toda $\tau$-functions $\tau_n(t,s),~n\in \BZ$ depend on two
sets of time-variables $t,s \in \BC^{\iy}$ and are defined by the
following bilinear identities, for all $m,n \in \BZ$:
\begin{equation}
\oint_{z=\iy}\tau_n(t-[z^{-1}],s)\tau_{m+1}(t'+[z^{-1}],s')
e^{\sum_1^{\iy}(t_i-t'_i)z^i}
z^{n-m-1}dz
\end{equation}
$$
=\oint_{z=0}\tau_{n+1}(t,s-[z])\tau_m(t',s'+[z])
e^{\sum_1^{\iy}(s_i-s'_i)z^{-i}}z^{n-m-1}dz,
$$
or, specified in terms of the Hirota symbol\footnote{for the
customary Hirota symbol $p(\pl_t)f\circ g:= p(\frac{\pl}{\pl
y})f(t+y)g(t-y)
\Bigl|_{y=0}$.}, by
\begin{eqnarray}
&&\sum_{j=0}^{\iy}p_{m-n+j}(-2a)p_j (\tilde\pl_t )
 e^{\sum_1^{\iy}(a_k\frac{\pl}{\pl t_k}+b_k \frac{\pl}{\pl s_k})}
 \tau_{m+1}\circ \tau_n \nonumber\\
&&=\sum_{j=0}^{\iy}p_{-m+n+j}(-2b)p_j (\tilde\pl_s )
 e^{\sum_1^{\iy}(a_k\frac{\pl}{\pl t_k}+b_k \frac{\pl}{\pl s_k})}
 \tau_{m}\circ \tau_{n+1}.
 \end{eqnarray}
%both, for the bi-infinite $(n,m\in\BZ)$
For the semi-infinite case, the same definitions hold, but for
$n,m\geq 0$.

\begin{theorem}
Two-Toda $\tau$-functions satisfy

\noindent (i) the {\bf KP-hierarchy} in $t$ and $s$ separately, of which
the first
equation reads:
$$
\left(\frac{\pl}{\pl t_1}\right)^4\log\tau
+6\left(\left(\frac{\pl}{\pl t_1}\right)^2\log\tau\right)^2
+3\left(\frac{\pl}{\pl t_2}\right)^2\log\tau
-4\frac{\pl^2}{\pl t_1\pl t_3}\log\tau = 0;
$$
\noindent (ii) an identity, involving $t,s$ and nearest
neighbors $\tau_{n-1}, \tau_n$ :
$$
\frac{\pl^2}{\pl s_1\pl
t_1}\log\tau_{n}=
-\frac{\tau_{n-1}\tau_{n+1}}{\tau_n^2};
$$
\noindent (iii) a (new) {\bf Identity} involving $t,s$ and nearest neighbors
$\tau_{n-1}, \tau_n$:
\be
\frac{\pl^2}{\pl s_2\pl
t_1}\log\tau_{n}=-2\frac{\pl}{\pl
s_1}\log\frac{\tau_{n}}{\tau_{n-1}}~.\frac{\pl^{2}}{\pl s_{1}\pl
t_{1}}\log\tau_{n}- \frac{\pl^{3}}{\pl s_{1}^2\pl
t_{1}}\log\tau_{n}.
\ee

\end{theorem}

\noindent The proof of this theorem will be given later in this section.

 The ``wave vectors"
defined in terms\footnote{We have
\bean
\chi(z)&=&(\ldots,z^{-1},1,z^{1},\ldots)\mbox{\,\, in the bi-infinite
case}\\
&=&(1,z,z^2,\ldots)\mbox{\,\, in the semi-infinite case.}
\eean
Also $\Lb^{-1}$ should always be interpreted as $\Lb^{\top}$ in the
semi-infinite case.}of $\tau_n(t,s)$,
\bea
\Psi_1(t,s,z)&=&\biggl(
        {\tau_n(t-[z^{-1}],s)\over\tau_n(t,s)}
        e^{\sum^{\iy}_1 t_iz^i}z^n
\biggr)_{n\in\BZ}=:e^{\sum^{\iy}_{1} t_i z^i}S_1 \chi(z) \nonumber\\
 \Psi^*_2(t,s,z)&=&\Biggl(\frac{ \tau_n(t,s+[z])}
{\tau_{n+1}(t,s)}e^{-\sum^{\iy}_1 s_iz^{-i}}
z^{-n}\Biggr)_{n\in\BZ}=:e^{-\sum ^{\iy}_{1} s_i z^{-i}}(S_2^{-1})^{\top}
\chi(z^{-1}),\nonumber\\
\eea
specify lower- and upper-triangular wave matrices $S_1$ and $S_2$
respectively. They, in turn,
 define a pair of matrices $L_1$ and $L_2$
\be
L_1:=S_1\Lambda S_1^{-1}=
 \sum_{-\iy<i\le0}a^{(1)}_i\Lb^i+\Lb
 ,~~L_2:=S_2\Lambda^{\top} S_2^{-1}
=\sum_{-1\le i<\iy}a^{(2)}_i\Lb^i ,
\ee
where $\Lb=(\dt_{j-i,1})_{i,j\in\BZ}$, and $a_i^{(1)}$ and
$a_i^{(2)}$ are diagonal matrices depending on $t=(t_1,t_2,\dots)$
and $s=(s_1,s_2,\dots)$. Then
\be
z\Psi_1=L_1 \Psi_1   ~~~\mbox{and}~~~ z^{-1}\Psi^*_2=L_2^{\top}
\Psi_2^*,
\ee
and the matrices $L_i$ satisfy the 2-Toda lattice equations:
 \be \frac{\pl
L_i}{\pl t_n}=\bigl[\bigl(L_1^n\bigr)_+,L_i\bigr]
\quad\hbox{and}\quad \frac{\pl L_i}{\pl
s_n}=\bigl[\bigl(L_2^n\bigr)_-,L_i\bigr] \quad
i=1,2~\mbox{and}~n=1,2,\dots
\ee
 %\hfill ({\bf Two-Toda Lattice})

%\vspace{0.2cm}

\noindent with $\Psi_1$ and $\Psi_2^*$ satisfying the
differential equations:
 $$
 \frac{\pl \Psi_1}{\pl t_n}=(L_1^n)_+\Psi_1 ~~~~~~~~~~~~~~~
  \frac{\pl \Psi_1}{\pl s_n}=(L_2^n)_-\Psi_1
  $$
  $$
 \frac{\pl \Psi^*_2}{\pl t_n}=-((L_1^n)_+)^{\top}\Psi^*_2
 ~~~~~~~~~~~~~~~
  \frac{\pl \Psi^*_2}{\pl
  s_n}=-((L_2^n)_-)^{\top}\Psi_2^*.
  $$

\vspace{0.2cm}

\noindent For future use, define the  diagonal
matrix:
\be
h:=(...,h_{-1},h_0,h_1,...), \mbox{ where } h_k(t,s):=
\frac{\tau_{k+1}(t,s)}{\tau_k(t,s)}. \ee In
\cite{AvM4}, we have shown that $L_1^k$ has the
following expression in terms of
$\tau$-functions\footnote{ $p_{\ell}(\tilde \pl)f
\circ g$ refers to the  Hirota operation, defined
before. Here the $p_{\ell}$ are the elementary Schur
polynomials $e^{\sum^{\iy}_{1}t_iz^i}:=\sum_{i\geq 0}
p_i(t)z^i$. Also $p_{\ell}(\tilde
\pl_t):=p_{\ell}(\frac{\pl}{\pl
t_1},\frac{1}{2}\frac{\pl}{\pl
t_2},\frac{1}{3}\frac{\pl}{\pl t_3},\ldots)$ and
 $p_{\ell}(-\tilde
\pl_s):=p_{\ell}(-\frac{\pl}{\pl
s_1},-\frac{1}{2}\frac{\pl}{\pl
s_2},-\frac{1}{3}\frac{\pl}{\pl s_3},\ldots).$}, \bea
L_1^k&=&\sum_{\ell=0}^{\iy}\mbox{diag}~
 \left(\frac{p_{\ell}(\tilde\pl_t)
\tau_{n+k-\ell+1}\circ\tau_n}
{\tau_{n+k-\ell+1} \tau_n}\right)_{n \in \BZ}\Lb^{k-\ell}
\nonumber\\
hL_2^{\top k}h^{-1}&=&
\sum_{\ell=0}^{\iy}\mbox{diag}~
 \left(\frac{p_{\ell}(-\tilde\pl_s)
\tau_{n+k-\ell+1}\circ\tau_n}
{\tau_{n+k-\ell+1} \tau_n}\right)_{n\in\BZ}\Lb^{k-\ell}.
\nonumber\\
 \eea

There is a general involution in the equation, which we shall
frequently use, namely $t\longleftrightarrow -s$, $L_1\longleftrightarrow
 hL_2^{\top}h^{-1}$.

 Finally, we
define the 2-Toda vertex operator, which is the generating function
for the algebra of symmetries, acting on $\tau$-functions (it will
play a role later!):
\be
\BX_{12}(t,s;u,v)=\Lb^{-1}e^{\sum_1^{\iy}(t_i u^i-s_i v^{i})}
e^{-\sum_1^{\iy}(\frac{ u^{-i}}{i}
\frac{\pl}{\pl t_i}-
 \frac{ v^{-i}}{i}
\frac{\pl}{\pl s_i})}\chi(uv),
\ee
 leading to a Virasoro algebra (see Appendix 1) with $\beta =1$ and thus
with   central
 charge $c=-2$,
\begin{eqnarray}
\frac{\pl}{\pl u}u^{k+1}\BX_{12}(t,s;u,v)&=&
 \left[\left.{}^{\beta}\BJ_k^{(2)}(t)\right|_{\beta=1},
\BX_{12}(t,s;u,v)\right]~,\nonumber\\
u^{k}\BX_{12}(t,s;u,v)&=&
 \left[\left.{}^{\beta}\BJ_k^{(1)}(t)\right|_{\beta=1},
\BX_{12}(t,s;u,v)\right]~,
\end{eqnarray}
 with generators (in $t$), explicitly given by (5.0.3).
Similarly, the involution $u\leftrightarrow v$,
 $t\leftrightarrow -s$ leads to the same Virasoro algebra in $s$, with
 same central charge.

\begin{proposition}

\begin{eqnarray}
\left(L_1^k\right)_{n,n}&=&\frac{p_{k}
 (\tilde\pl_t )\tau_{n+1}
 \circ \tau_n}{\tau_{n+1} \tau_n}
 ~=~~\frac{\pl}{\pl t_k} \log \frac{\tau_{n+1}}{\tau_{n}}
\nonumber\\
\left(hL_2^{\top k}h^{-1}\right)_{n,n}
& =&\frac{p_{k}
 (-\tilde\pl_s )\tau_{n+1}
 \circ \tau_n}{\tau_{n+1} \tau_n}
 ~=~~-\frac{\pl}{\pl s_k} \log \frac{\tau_{n+1}}{\tau_{n}}
\end{eqnarray}
 and
\begin{eqnarray}
\left(L_1^k\right)_{n,n+1}&=&\frac{p_{k-1}
 (\tilde\pl_t )\tau_{n+2}
 \circ \tau_n}{\tau_{n+2} \tau_n}
 ~=~~\frac{\frac{\pl^2 \log \tau_{n+1}}{\pl s_1 \pl t_{k}}}
 {\frac{\pl^2 \log \tau_{n+1}}{\pl s_1 \pl t_{1}}}
\nonumber\\
\left(hL_2^{\top k}h^{-1}\right)_{n,n+1}
& =&\frac{p_{k-1}(-\tilde\pl_s )\tau_{n+2}
 \circ \tau_n}{\tau_{n+2} \tau_n}
 =~\frac{\frac{\pl^2 \log \tau_{n+1}}{\pl t_1 \pl s_{k}}}
 {\frac{\pl^2 \log \tau_{n+1}}{\pl s_1 \pl t_{1}}}.
\end{eqnarray}
\end{proposition}

\proof Relations (2.1.12) follows from (2.1.7) and a
 standard argument; see \cite{AvM4}, Theorem 0.1, formula (0.15).

To prove (2.1.13), set $m=n+1$, all $b_k$ and $a_k=0$, except for
one
 $a_{j+1}$, in the Hirota bilinear relation (2.1.2). The first nonzero term
in the sum on the left hand
 side of that relation, which is also the only one containing
 $a_{j+1}$ linearly, reads
 \begin{equation}
 p_{j+1}(-2a)p_j(\tilde\pl_t)
 e^{a_{j+1}\frac{\pl}{\pl t_{j+1}}}\tau_{n+2}\circ \tau_n+...
  =-2a_{j+1}p_j(\tilde\pl_t)\tau_{n+2}\circ
  \tau_n+O(a_{j+1}^2),
  \end{equation}
  whereas the right hand side equals
  \begin{equation}
  p_{0}(0)p_1(\tilde\pl_s)
 e^{a_{j+1}\frac{\pl}{\pl t_{j+1}}}\tau_{n+1}\circ \tau_{n+1}
 =\frac{\pl}{\pl s_1}(1+a_{j+1}\frac{\pl}{\pl t_{j+1}}+...)
 \tau_{n+1}\circ \tau_{n+1}.
 \end{equation}
Comparing the coefficients of $a_{j+1}$ in (2.1.14) and (2.1.15)
yields
$$
-2~p_j(\tilde\pl_t)\tau_{n+2}\circ \tau_n=
\frac{\pl^2}{\pl s_1 \pl t_{j+1}} \tau_{n+1}\circ \tau_{n+1}
;$$
in particular, we have
\begin{equation}
\frac{p_{k-1}(\tilde\pl_t)\tau_{n+2}\circ \tau_n}
{\tau_{n+1}^2}=-\frac{\pl^2}{\pl s_1 \pl t_k}\log \tau_{n+1} ,
 \end{equation}
and so, for $k=1$,
\begin{equation}
\frac{\tau_n \tau_{n+2} }
{\tau_{n+1}^2}=-\frac{\pl^2}{\pl s_1 \pl t_1}\log \tau_{n+1} .
 \end{equation}
Dividing (2.1.16) and (2.1.17) leads to the
 first equality in (2.1.13), since according to (2.1.9), the
 $(n,n+1)$-entry of $L_1^k$ is precisely given by (2.1.16).
The similar result for $L_2^k$ is given by the involution
 $$
 t \longleftrightarrow  -s ~~\mbox{and}~~  L_1\longleftrightarrow
 hL_2^{\top }h^{-1}.
 $$
\qed

\begin{lemma}
The first upper-subdiagonal of $L^2_{1}$ and
$hL_{2}^{\top 2}h^{-1}$ reads:

\begin{eqnarray}
\left(L_1^2\right)_{n,n+1}&=&\frac{\pl}
 {\pl t_1}\log\frac{\tau_{n+2}}{ \tau_n}\nonumber\\
&=&\frac{\frac{\pl^2} {\pl s_1\pl
t_2}\log\tau_{n+1}}{\frac{\pl^2}{\pl s_1\pl t_1}
 \log\tau_{n+1}}
 \nonumber\\&=&
\frac{\pl}{\pl t_1}\log
 \left( \left(\frac{\tau_{n+1}}{\tau_n } \right)^2
\frac{\pl^2}{\pl s_1\pl
t_1}\log\tau_{n+1}  \right)  \nonumber\\
 && \nonumber\\
\left(hL_2^{\top 2}h^{-1}\right)_{n,n+1}
& =&
 -\frac{\pl}
 {\pl s_1}\log\frac{\tau_{n+2}}{ \tau_n}\nonumber\\
 & =&\frac{\frac{\pl^2} {\pl
t_1\pl s_2}\log\tau_{n+1}}{\frac{\pl^2}{\pl t_1\pl s_1}
 \log\tau_{n+1}} \nonumber\\
 &=&
-\frac{\pl}{\pl s_1}\log
 \left( \left(\frac{\tau_{n+1}}{\tau_n } \right)^2
\frac{\pl^2}{\pl s_1\pl
t_1}\log\tau_{n+1}  \right) .
\end{eqnarray}
\end{lemma}

\proof From Proposition 2.2 $(k=2)$, we have the first
two identities in (2.1.18); it also follows from these
identities that
\begin{eqnarray*}
\frac{\pl^2\log\tau_{n+1}}{\pl s_1\pl
t_2}
 %&=&\frac{\pl^2\log\tau_{n+1}}{\pl s_1\pl t_1}
% \frac{\frac{\pl}{\pl t_1}(\tau_{n+2}\circ\tau_n)}{\tau_{n+2}\tau_n}\\
&=&\frac{\pl^2\log\tau_{n+1}}{\pl s_1\pl
t_1}\frac{\pl}{\pl t_1}\log\frac{\tau_{n+2}}{\tau_n}\\
&=&\frac{\pl^2\log\tau_{n+1}}{\pl s_1\pl
t_1}\left(\frac{\pl}{\pl t_1}\log
\frac{\tau_{n+2}}{\tau_{n+1}}+\frac{\pl}{\pl t_1}\log
\frac{\tau_{n+1}}{\tau_n}\right)\\
&=&\frac{\pl^2\log\tau_{n+1}}{\pl s_1\pl
t_1}\left(\frac{\pl}{\pl
t_1}\log\left(-\frac{\tau_{n+1}}{\tau_n}\frac{\pl^2}{\pl
s_1\pl t_1}\log\tau_{n+1}\right)\right.\\ &
&\hspace{1cm}+\,\left.\frac{\pl}{\pl t_1}\log\frac{
\tau_{n+1}}{\tau_n}\right),\mbox{using (2.1.17)}\\
&=&\frac{\pl^2\log\tau_{n+1}}{\pl s_1\pl
t_1}\left(2\frac{\pl}{\pl t_1}\log\frac{
\tau_{n+1}}{\tau_n}+\frac{\pl}{\pl
t_1}\log\left(\frac{\pl^{2}}{\pl s_{1}\pl
t_{1}}\log\tau_{n+1}\right)\right)\\
&=&2\frac{\pl}{\pl t_1}\log\frac{\tau_{n+1}}{\tau_n}
 \frac{\pl^{2}}{\pl s_{1}\pl t_{1}}\log\tau_{n+1}
 +\frac{\pl}{\pl t_1}\left(\frac{\pl^{2}}{\pl
s_{1}\pl t_{1}}\log\tau_{n+1}\right),
\end{eqnarray*}
which establishes the first equation (2.1.18). The second
 equation (2.1.18) is simply the dual of the first one by
 $t_i\mapsto -s_i$.
\qed

\noindent\underline{\sl Proof of Theorem 2.1}:  The first
 statement concerning the KP hierarchy is standard.
The proof of the second identity
 follows immediately from (2.1.16), for $k=1$, and the
 third identity from the
 last identity  in the proof of Lemma 2.3 and the duality.\qed

\medbreak
A {\em prominent example} of the semi-infinite 2-Toda lattice is given by
an (arbitrary) $(t,s)$-dependent  semi-infinite matrix
\be
m_{\iy}(t,s)=\left( \mu_{ij}(t,s)\right)_{0\leq i,j
<\iy },~\mbox{with}~m_{n}(t,s)=\left(
\mu_{ij}(t,s)\right)_{0\leq i,j\leq n-1},\ee
 evolving according to the equations
 \be
 \frac{\pl m_{\iy}}{\pl t_k}=\Lb^k m_{\iy}~~\mbox{and}~~
\frac{\pl m_{\iy}}{\pl s_k}=- m_{\iy} (\Lb^{\top})^k.
\ee
 According to \cite{AvM2}, the formal solution to this problem is given by
\be
m_{\iy}(t,s)=e^{\sum_1^{\iy} t_i \Lb^i}m_{\iy}(0,0)
 e^{-\sum_1^{\iy} s_i \Lb^{\top i}}=S_1^{-1}(t,s)S_2(t,s),
\ee where the associated unique factorization into
lower- times upper-triangular matrices actually lead
to the wave matrices $S_1$ and $S_2$, as defined in
the general Toda theory (2.1.4). The expression
(2.1.21) contains the matrix of Schur
polynomials\footnote{The Schur polynomials $p_i$,
defined by $e^{\sum^{\infty}_1 t_i z^i
}=\sum^{\infty}_0 p_k(t)z^k$ and $p_k(t)=0$ for $k<0$,
are not to be confused with the bi-orthogonal
polynomials $p_i^{(k)}$, $k=1,2$.}

$$
 e^{\sum_{1}^{\iy}t_i \Lb^i}=\sum_0^{\iy} \Lb^i p_i(t)=
\Bigl (p_{j-i}(t)\Bigr)_{\tiny\begin{array}{l}
\scriptsize1\leq i < \iy\\ \scriptsize1\leq j<\iy
\end{array}},
 $$
of which a truncated version is given by the following
$n\times\iy$ submatrix $p_i(t)$:
\begin{eqnarray}
E_n(t)&=&\left(\begin{array}{ccccc|cc}
1&p_1(t)&p_2(t)&\pp&p_{n-1}(t)&p_n(t)&\pp\\
0&1&p_1(t)&\pp&p_{n-2}(t)&p_{n-1}(t)&\pp\\
\vdots&\vdots& & & \vdots&\vdots& \\
0&0&0&\pp&p_1(t)&p_2(t)&\pp\\ 0&0&0&\pp&1&p_1(t)&\pp
\end{array}\right)
\nonumber\\ &=&\Bigl(p_{j-i}(t)\Bigr)_{ 1\leq i\leq n
\atop 1\leq j<\iy }
\end{eqnarray}
So, for a semi-infinite initial condition
$m_{\iy}(0,0)$, the $\tau$-functions of the 2-Toda
problem are given by
\be
\tau_n(t,s):=\det m_{n}(t,s)
 =\det \left( E_n(t)~m_{\iy}(0,0)~E_n^{\top}(-s)\right)
.\ee Incidentally, the wave vectors $\Psi_1$ and
$\Psi_2^*$ define monic polynomials $p^{(1)}(x)$ and
$p^{(2)}(y)$,
\be
\begin{tabular}{lll}
$\Psi_1:=e^{\sum t_kz^k}p^{(1)}(z)$&\mbox{and}&
 $\Psi_2^*:=e^{-\sum
s_kz^{-k}}h^{-1}p^{(2)}(z^{-1}) $   \\ $\,\,\,~~~=e^{\sum
t_kz^k}S_1\chi(z)$& &$\,\,\,~~~=e^{-\sum
s_kz^{-k}}(S_2^{-1})^{\top}\chi(z^{-1}),$
\end{tabular}
\ee
which are {\em
bi-orthogonal} with regard to the original matrix
$m_{\iy}$; that is, for all $t,s$:
 \be\la
p_n^{(1)},p_m^{(2)}\ra
=\dt_{n,m}h_n~~\mbox{for the inner-product defined by}~~
 \la x^i,y^j  \ra:= \mu_{ij},
 \ee
with $h_n$ as in (2.1.8).

%\newpage
\vspace{.5cm}

\noindent{\bf Proposition } {\it Given the
semi-infinite initial condition $m_{\iy}(0,0)$, the
$2$-Toda $\tau$-function has the following expansion
in Schur polynomials\footnote{For a given Young
diagram $\lb_1 \geq ...\geq\lb_n$, define
$s_{\lb}(t)=\det (p_{\lb_i-i+j}(t))_{1\leq i,j \leq
n}$. },
 \be
\tau_n(t,s)
% =\det \left(
%E_n(t)~m~E_n^{\top}(-s)\right)
 = \sum_{\lb,~ \nu \atop
\hat \lb_1,~\hat \nu_1 \leq n}
    \det(m_{}^{\lb,\nu}) s_{\lb}(t) s_{\nu}(-s),
 ~~ \mbox{for $n>0$}, \ee
where the sum is taken over all Young diagrams $\lb$
and $\nu$, with first columns $\hat\lb_1$ and
$\hat\nu_1\leq n $ and where
  \be
   m_{}^{\lb,\nu}:=\left( \mu_{\lb_{i}-i+n,\nu_{j}-j+n}
    \right)_{1\leq i,j\leq n}.
 \ee
 }

\proof Note that every increasing sequence $1\leq
k_1<...<k_n<\iy$ can be mapped into a Young diagram
$\lb_1\geq \lb_2\geq...\geq \lb_n\geq 0$, by setting
$k_j=j+\lb_{n+1-j}$. Relabeling the indices $i,j$ with
$1\leq i,j\leq n$, by setting $j':=n-j+1,~i':=n-i+1$,
we have $1\leq i',j'\leq n$ and $k_j-i=\lb_{j'}-j'+i'$
and $k_i-1=\lb_{i'}-i'+n$. The same can be done for
the sequence $1\leq \ell_1<...<\ell_n<\iy$, leading to
the Young diagram $\nu$, using the same relabeling.
 Applying the Cauchy-Binet formula twice, expression (2.1.23) leads to:
  \bean
 \lefteqn{\tau_n(t,s)}\\ %&&\\
 &=& \det \left( E_n(t)~m_{\iy}(0,0)
  ~E_n^{\top}(-s)\right)\\
 &=& %\hspace{-.3cm}
    \sum_{1\leq k_1<...<k_n<\iy}\det
    \left(p_{k_j-i}(t)\right)_{1\leq i,j \leq n}
    \det
    \left((m_{\iy}(0,0)E_n^{\top}(-s))_{k_i,\ell}
    \right)_{1\leq i,\ell\leq n}\\
 &=& \hspace{-.3cm}
    \sum_{1\leq k_1<...<k_n<\iy}\det
    \left(p_{k_j-i}(t)\right)_{1\leq i,j \leq n}
    \det
    \left( (\mu_{k_i-1,j-1})_{ 1\leq i\leq n\atop
    1\leq j<\iy }
      (p_{i-j}(-s))_{ 1\leq i<\iy \atop 1\leq j\leq n
     })\right) \\
 &=& %\hspace{-.3cm}
   \sum_{1\leq k_1<...<k_n<\iy}\det
    \left(p_{k_j-i}(t)\right)_{1\leq i,j \leq n}\\
 &&~~~~~~  \sum_{1\leq \ell_1<...<\ell_n<\iy}\det
    \left( \mu_{k_i-1,\ell_j-1}\right)_{1\leq i,j \leq n}
    \det
    \left( p_{\ell_i-j}(-s)\right)_{1\leq i,j \leq n}
    \\
 &=&  %\hspace{-.3cm}
   \sum_{\lb \atop \hat \lb_1\leq n}\det
    \left(p_{\lb_{j'}-j^{\prime}+i'}(t)
    \right)_{1\leq i',j^{\prime} \leq n}
    \\
 &&~~~~~~  \sum_{\nu \atop \hat \nu_1\leq n}\det
    \left( \mu_{\lb_{i'}-i'+n,\nu_{j'}-j'+n}
    \right)_{1\leq i',j' \leq n}
    \det
    \left( p_{\nu_{i'}-i^{\prime}+j'}(-s)
    \right)_{1\leq i',j' \leq n}
    \\
 &=&
    \sum_{\lb,~ \nu \atop \hat \lb_1,~\hat \nu_1 \leq n}
    \det(\mu_{}^{\lb,\nu}) s_{\lb}(t) s_{\nu}(-s).
\eean

\vspace{-1cm}\qed

%\newpage

\subsection{Reduction to H\"ankel matrices:
the standard Toda lattice and a Virasoro algebra
 of constraints}

 In the
 notation of (2.1.7), consider the locus of $(L_1,L_2)$'s
 such that $L_1=L_2$. This means the matrix $L_1=L_2$ is
 tridiagonal. From the equations (2.1.7), it follows
 that along that locus, $$ \frac{\pl(L_1-L_2)}{\pl
t_n}=0\quad\quad\quad\frac{\pl(L_1-L_2)}{\pl s_n}=0. $$
 We now
define new variables $t'_n$ and $s'_n$ by
\be
t'_n=t_n-s_n ~~\mbox{and}~~s'_n=t_n+s_n,
\ee
and thus
 \be
\frac{\pl}{\pl t'_n}=\frac{1}{2}\left(\frac{\pl}{\pl t_n}-\frac{\pl}{\pl
s_n}\right),\quad\quad \frac{\pl}{\pl
s'_n}=\frac{1}{2}\left(\frac{\pl}{\pl t_n}+\frac{\pl}{\pl
s_n}\right).
 \ee
Then
\begin{equation}
\frac{\pl L_i}{\pl s'_n}=\frac{1}{2}\left(\frac{\pl}{\pl t_n}+\frac{\pl}{\pl
s_n}\right)L_i=[(L_1^n)_++(L_2^n)_{-},L_i]=[L_i^n,L_i]=0
 \end{equation}
and
\begin{eqnarray}
\frac{\pl L_i}{\pl t'_n}&=&\frac{1}{2}\left(\frac{\pl}{\pl t_n}-\frac{\pl}{\pl
s_n}\right)L_i\nonumber\\
 &=&\frac{1}{2}[(L_1^n)_+-(L_2^n)_{-}+L_i^n,L_i] \nonumber\\
 &=&
 [(L_1^n)_+,L_i], ~\mbox{using}~L_1=L_2.
\end{eqnarray}
So, equation (2.2.3) implies that $L_1=L_2$ is independent of $s'$.
Since $\tau(t,s)$ is a function of $t-s$ only, we may set
$\tau(t'):=\tau(t-s)$. %Define
%$h:=~\mbox{diag}(...h_{-1},h_0,h_1,...)$, with
%$h_n(t):=\tau_{n+1}(t)/\tau_n(t)$.
After noting (see (2.1.12))
$$
\frac{\pl}{\pl t'_k}\log h_{n}  = \left(L_1^k\right)_{n,n},$$
 this situation leads to the standard Toda lattice equations on
  symmetric and tridiagonal matrices
 $L=h^{-1/2}L_1h^{1/2}$ and wave vectors (expressed in
 terms of the 2-Toda wave vectors (2.1.4))
 \begin{eqnarray}
  \Psi(t',z) :=h^{-1/2}\Psi_1 (t,s;z)
e^{-\frac{1}{2}\sum^{\iy}_{1}(t_i+s_i)z^i}&=&h^{1/2}\Psi_2^\ast
(t,s;z^{-1})e^{\frac{1}{2} \sum^{\iy}_{1}(t_i+s_i)z^i};\nonumber\\
&=&e^{\frac{1}{2}\Sigma
t'_iz^i}\left(z^n\frac{\tau_n(t'-[z^{-1}])}{\sqrt{\tau_n(t')\tau_{n+1}(t')}}
\right)_{n\geq 0}\nonumber\\
& &
\end{eqnarray}
 namely
\be
 L\Psi=z\Psi ,~~
 \frac{\partial \Psi}{\partial t'_n} =
\frac{1}{2} (L^n)_{sk}\Psi~~\mbox {and }~~
 \frac{\partial L}{\partial t'_n} =
\frac{1}{2}\left[ (L^n)_{sk}, L\right],
 %\quad\quad\mbox{with}
 %\quad L = h^{-\frac{1}{2}} \, L_1
%h^{\frac{1}{2}}
\ee
\hfill  {\bf (Standard Toda lattice)}

 \noindent where $()_{sk}$ refers to the
 skew-part in the skew and lower-triangular
 Lie decomposition.

We now define the standard Toda vertex operator as the
reduction of the two-Toda vertex operator $\BX_{12}$,
defined in (2.1.10), using (2.2.1) and (2.2.2):
\be
\BX(t';z):= \BX_{12}(t,s;z,z) \Bigr|_{
{t=\frac{s'+t'}{2}}\atop{s=\frac{s'-t'}{2} }  }
=\Lambda^{\top}\chi(z^2)e^{\sum_1^{\infty}
t'_iz^i}e^{-2\sum_1^{\infty} \frac{z^{-i}}{i}
\frac{\partial}{\partial t'_i}}. \ee In the rest of
this section, we shall omit ${}^{\prime}$ in $t'$ and
$s'$. This vertex operator $\BX(t;z)$ generates a
Virasoro algebra, with $\beta =2$ and thus central
charge $c=1$ (see Appendix 1 and (5.0.5))
\be
\frac{d}{du}u^{k+1}\BX(t,u)
 =\left[ \left.{}^{\beta}\BJ^{(2)}_k(t)\right|_{\beta=2},\BX(t;u)\right],
~~k \in \BZ \, ;
\ee

\vspace{0.6cm}

An interesting {\em semi-infinite example} of the
standard Toda lattice is obtained by considering a
weight $\rho(z)dz=e^{-V(z)}dz$ defined on an interval
$F\subset \BR$, satisfying (0.0.6) and an
inner-product
\be
\la f,g\ra_t:=\int_{\BR} f(z)g(z) e^{\sum_1^{\iy}t_iz^i}
 \rho(z) dz,
\ee
leading to a $t$-dependent moment matrix
\be
m_{\iy}(t)=\left(\mu_{ij}(t)\right)_{0\leq i,j\leq \iy}
 =\left(\la y^i,y^j \ra_t\right)_{0\leq i,j < \iy}.
\ee
 \hfill ({\bf H\"ankel matrix})

\begin{theorem}{\em (Adler-van Moerbeke \cite{AvM1})}\,
The vector $\tau(t)=(\tau_0=1,\tau_1(t),...)$ of integrals \bea
\tau_n(t)&:=&\frac{1}{n!} \int_{{\cal H}_n}
 e^{\Tr (-V(M)+\sum_1^{\iy}t_i M^i)}dM\nonumber\\
 &=& \frac{1}{n!} \int_{F^n} \Delta_n(z)^2\prod_{k=1}^n
  e^{\sum_{i=1}^{\iy} t_i z_k^i} \rho(z_k) dz_k\nonumber\\
  &=& \det \left(\mu_{ij}(t)\right)_{0\leq i,j\leq n-1}
\eea is a set of $\tau$-functions for the standard Toda lattice.
Also, each $\tau_n(t)$ satisfies the KP hierarchy, of which the
first equation is given in Theorem 2.1. It
also satisfies (i) of theorem 0.2, with $t_1=x$. Moreover, the $\tau_n$'s
satisfy the following Virasoro constraints for $\beta=2$,(see
(2.2.8)),
  \bea
  0&=&\JR^{(2)}_{m}\tau(t)
 ,\quad
m\geq -1,\nonumber\\ &:=&\sum_{k
\geq 0}\left.
\left(
 -a_k
 \sum_{i+j=k+m}:  {}^{\beta}\BJ^{(1)}_{i}  {}^{\beta}\BJ^{(1)}_{j}:
+b_k ~\,^{\beta}\BJ^{(1)}_{k+m+1} \right)\right|_{\beta
=2}\tau(t)\nonumber\\ &=& \left( \sum_{k\geq 0}\Bigl(-a_k
(~{}^{\beta}J_{k+m}^{(2)} +2n\,\, ^{\beta}J_{k+m}^{(1)} +n^2
J_{k+m}^{(0)}
 )\right.\nonumber\\
 &&~~~~~~~~~~~~~~\left. +b_k ( \, ^{\beta}J_{k+m+1}^{(1)}  +
 n \dt_{k+m+1,0}   )\Bigr)\Bigl|_{\beta=2}\tau_n(t)\right)_{n \geq
0}\nonumber\\ \eea where the $a_{k}$ and $b_{k}$ are
the coefficients (0.0.6) of the rational function
$\rho'/\rho$, and where  the $^{\beta}\BJ^{(i)}_k$ and
$^{\beta}J^{(i)}_k$ for $\beta =2$, are given by
 (5.0.1), (5.0.3) and (5.0.5). The
relation between the vertex operator $\BX(t,u)$ and ${\cal
J}_{m}^{(2)}$ is given by
%in terms of (2.2.9) and the vertex operator $\BX(t,u)$ in
%(2.2.6), we have for $m\geq -1$,
 \be
 \frac{\pl}{\pl u}u^{m+1}f(u)\BX_{}(t,u)\rho(u)
 =\left[ -{\cal
J}^{(2)}_{m},\BX_{}(t,u)\rho(u)\right]. \ee

\end{theorem}

\medskip\noindent{\it Sketch of proof:\/} Formula (2.2.13) is a
direct consequence of (2.2.8), while (2.2.12)
hinges on the fact that the vector $I=(\tau_0,\tau_1,...,n!\tau_n,...)$
of Toda lattice $\tau$-functions is a fixed point for a certain
integrated vertex operator, in the following sense: $$
\left({\cal Y}I\right)_n:=I_n,~\mbox{for}~n\geq 1,~~\mbox{where}~
 {\cal Y}:=\int_F  du \rho(u) \BX(t,u),
 $$
which is just an iterated integral formula. Upon
integrating (2.2.13) on the full range $F$, deduce
$[\JR_m^{(2)}, {\cal Y} ] =  0 $ from the boundary
conditions (0.0.6). Acting on $I$ with this relation,
one deduces (2.2.12) by induction on $n$ and the fact
that $\tau_0=1$ . \qed

\subsection{Reduction to Toeplitz matrices: two-Toda
 Lattice and an
SL(2,$\BZ$)-algebra of constraints }

Consider the following inner product, depending on $(t,s)$,
\be
 \la f(z),g(z)\ra_{t,s}:=\oint_{S^{1}} \frac{\rho(dz)}{2\pi i z}f(z)g(z^{-1})
 e^{\sum_1^{\iy}(t_iz^i-s_iz^{-i})},
 \ee
where the integral is taken over the unit circle $S^1$ around the
origin in the complex plane $\BC$. Instead of having
$z^{k\top}=z^k$ in the H\"ankel inner-product, we have
$z^{k\top}=z^{-k}$ in this inner-product:
\be
\la z^kf(z),g(z)\ra_{t,s} = \la
f(z),z^{-k}g(z)\ra_{t,s}. \ee
 Thus the moment
matrix $m_{\iy}$, with entries
\be
\mu_{k\ell}(t,s)=\la z^k, z^{\ell}\ra_{t,s}=\oint_{S^1}
 \frac{\rho(z)dz}{2\pi i z}z^{k-\ell}e^{\sum_1^{\iy}(t_iz^i-s_iz^{-i})},
\ee is a Toeplitz matrix for all $t,s$, satisfying the
differential equations (2.1.20) of the 2-Toda lattice, i.e., $$
\frac{\pl \mu_{k\ell}}{\pl t_i}=\mu_{k+i,\ell} ~~\mbox{and}~~
\frac{\pl \mu_{k\ell}}{\pl s_i}=-\mu_{k,\ell +i}. $$

\begin{theorem}
For $\rho(dz)=dz$, the vector
 $\tau(t,s)=(\tau_0=1,\tau_1(t,s),...)$, with
\bea \tau_n(t,s)&=&\int_{U(n)}e^{\sum_1^{\iy}\Tr
(t_iM^i-s_i\bar M^i)}dM\nonumber\\ &=& \frac{1}{n!}
\int_{(S^1)^{n}}|\Dt_n(z)|^{2}
 \prod_{k=1}^n
\left(e^{\sum_1^{\iy}(t_i z_k^i-s_iz_k^{-i})}
 \frac{dz_k}{2\pi i z_k}\right)\nonumber\\&&\nonumber\\
&=&
\det \left( \mu_{k\ell }(t,s)\right)_{0\leq k,\ell \leq n-1}
 \nonumber  \\&&\nonumber\\ &=& \sum_{\{\mbox{\footnotesize
%partitions
Young diagrams $\lambda$ with first column $\leq n$
 }\}}s_{\lambda}(t)s_{\lambda}(-s)\nonumber\\ \eea
 is a vector of $\tau$-functions for the two-Toda lattice.
 Hence, they
 satisfy the identity (2.1.3), and (ii) of Theorem 0.2,
 with $t_1=x,~s_1=y$.
  Moreover, they are annihilated by
 the following algebra of
 three Virasoro partial differential operators, which form a
 $SL(2,\BZ)$- algebra, (note: the expression below is a
 semi-infinite vector),
 %in terms of ${\cal J}^{(2)}_k$ defined in ():
\bea 0&=&{\cal V}^{(2)}_k\tau(t,s)\quad
 \mbox{for}~
\left\{\begin{array}{l}
k =-1 ,~\theta= 0\\ k=0 , ~~\theta~~ {arbitrary}\\
  k=1 ,~
\theta= 1\end{array}\right\}
\mbox{\,\,only
 }  \nonumber\\
 &=& \left( \left.{}^{\beta}\BJ_{k}^{(2)}(t)-
  {}^{\beta} \BJ_{-k}^{(2)}(-s)
-k\bigl( \theta~~ {}^{\beta}\BJ_{k}^{(1)}(t)+(1-\theta)~~
 {}^{\beta}\BJ_{-k}^{(1)}(-s)
 \bigr)\right)\right|_{\beta=1}\tau(t,s)
, \nonumber\\ \eea where $\theta$ is an arbitrary parameter and
 $  {}^{\beta} \BJ_{k}^{(2)}$ for $\beta =1$
 is given in (5.0.4).
 This is a subalgebra of a Virasoro algebra
(generated by the 2-Toda vertex operator (2.1.10))
\be
\frac{d}{du}u^{k+1} \frac{\BX_{12}(t,s;u,u^{-1})}{u}
 =\left[ {\cal V}^{(2)}_k(t,s),
  \frac{\BX_{12}(t,s;u,u^{-1})}{u}\right],~~k \in \BZ,
\ee
 of central charge $c=0$:
\be
 \left[{\cal V}^{(2)}_k,{\cal V}^{(2)}_{\ell}\right]=
 (k-\ell){\cal V}^{(2)}_{k+\ell}.
 \ee

\end{theorem}

\vspace{0.4cm}

\noindent The statement hinges on the following statement about
 the vertex operator:

\begin{proposition}
For general weight $\rho(dz)=e^{V(z)}dz$, the column
vector of 2-Toda $\tau$-functions (slightly rescaled),
$\tau(t,s)=(\tau_0,\tau_1,...)~,~~\mbox{with}~~$
\be
I_n(t,s)=n!\tau_n =n!\int_{U(n)}e^{\sum_1^{\iy}\Tr
(t_iM^i-s_i\bar M^i)}e^{\Tr~V(M)}dM
  \ee
 is a fixed point, in the sense,
\be
\left({\cal
Y}_{}(t,s;\rho)I\right)_n=I_n,~\mbox{for}~n\geq 1,
\ee
for the operator
 \be
  {\cal Y}_{}(t,s;
\rho)=\int_{S^1} \frac{\rho(du)}{2\pi i u}
\BX_{12}(t,s;u,u^{-1}).
\ee

\end{proposition}

\vspace{0.6cm}

\noindent\underline{\sl Proof of Proposition 2.6 and Theorem 2.5}:
On the one hand, using the fact that $\bar z=1/z$ on the circle
$S^1$ and a property of Vandermonde determinants\footnote{The
following holds:
 $$\det(u_k^{\ell -1})_{1\leq\ell,k\leq N}
\det(v_k^{\ell -1})_{1\leq\ell,k\leq N}=\sum_{\sigma\in
S_N}\det\Bigl(u^{\ell-1}_{\sigma(k)}
v^{k-1}_{\sigma(k)}\Bigr)_{1\leq\ell,k\leq N}.$$} and Theorem 1.1, we have, for a
general weight $\rho(dz)=e^{V(z)}dz$,
 $\displaystyle{}$
\bea n!\tau_n(t,s)&=&n!\int_{U(n)}e^{\sum_1^{\iy}\Tr
(t_iM^i-s_i\bar M^i)}  e^{\Tr V(M)}dM\nonumber\\ &=&
 \int_{(S^1)^{n}}|\Dt_n(z)|^{2}
 \prod_{k=1}^n
\left(e^{\sum_1^{\iy}(t_i z_k^i-s_iz_k^{-i})}
 \frac{\rho(dz_k)}{2\pi i z_k}\right)\nonumber\\
&=&
 \int_{(S^1)^{n}}\Dt_n(z)\Dt_n(\bar z)
 \prod_{k=1}^n
\left(e^{\sum_1^{\iy}(t_i z_k^i-s_iz_k^{-i})}
 \frac{\rho(dz_k)}{2\pi i z_k}\right)\nonumber\\
&=&
 \int_{(S^1)^{n}}\sum_{\sigma\in S_n}
 \det\left(z_{\sigma(m)}^{\ell-1}
 \bar z_{\sigma(m)}^{m-1}  \right)_{1\leq \ell,m\leq n}
 \prod_{k=1}^n
\left(e^{\sum_1^{\iy}(t_i z_k^i-s_iz_k^{-i})}
 \frac{\rho(dz_k)}{2\pi i z_k}\right)
 \nonumber\\
&=&
 \sum_{\sigma\in S_n} \det\left(\int_{S^1}z_{}^{\ell-1}
 \bar z_{}^{m-1}
e^{\sum_1^{\iy}(t_i z^i-s_iz^{-i})}
 \frac{\rho(dz)}{2\pi i z}\right)_{1\leq \ell,m\leq n}
  \nonumber\\
  &=&
 n!\det\left(\int_{S^1}z^{\ell-m}
e^{\sum_1^{\iy}(t_i z^i-s_iz^{-i})}
 \frac{\rho(dz)}{2\pi i z}\right)_{1\leq \ell,m\leq n}\nonumber\\
&=& n!\det \left( \mu_{\ell m}(t,s)\right)_{0\leq
\ell,m \leq n-1}, \eea
 yielding the third equality of
(2.3.4). On the other hand, for $n\geq 1$, we have
 \bea I_n(t,s)&=&n!\tau_n(t,s)\nonumber\\
&=&\int_{(S^1)^{n}}|\Dt_n(z)|^{2}
 \prod_{k=1}^n
\left(e^{\sum_1^{\iy}(t_i z_k^i-s_iz_k^{-i})}
 \frac{\rho(dz_k)}{2\pi i z_k}\right)
\nonumber\\
&&\nonumber\\
 &=&\int_{S^1}\frac{\rho(du)}{2\pi iu}e^{\sum_1^{\iy}(t_i u^i-s_i u^{-i})}
 u^{n-1}u^{-n+1}\nonumber\\
 &&\int_{(S^1)^{n-1}} \Dt_{n-1}(z)\bar \Dt_{n-1}(z)
 \prod_{k=1}^{n-1}\left(1-\frac{z_k}{u}\right)
  \left(1-\frac{u}{z_k}\right)
e^{\sum_1^{\iy}(t_i z_k^i-s_i z_k^{-i})}
\frac{\rho(dz_k)}{2\pi i z_k}\nonumber\\
&=&\int_{S^1}\frac{\rho(du)}{2\pi iu}e^{\sum_1^{\iy}(t_i u^i-s_i
u^{-i})} ~e^{- \sum_1^{\iy}\left(\frac{u^{-i}}{i}
 \frac{\pl}{\pl t_i}
 -\frac{u^{i}}{i}\frac{\pl}{\pl s_i}\right)}
 \nonumber\\
&&~~~~~\int_{(S^1)^{n-1}}
  \Dt_{n-1}(z)\bar \Dt_{n-1}(z)
 \prod_{k=1}^{n-1}
e^{\sum_1^{\iy}(t_i z_k^i-s_i z_k^{-i})}
\frac{\rho(dz_k)}{2\pi iz_k}\nonumber\\
&=&\int_{S^1}\frac{\rho(du)}{2\pi iu}e^{\sum_1^{\iy}(t_i u^i-s_i
u^{-i})} e^{- \sum_1^{\iy}\left(\frac{u^{-i}}{i}\frac{\pl}{\pl t_i}
 -\frac{u^{i}}{i}\frac{\pl}{\pl s_i}\right)}
 I_{n-1}(t,s)
%&=&\int_{\BR}du\rho(u)I_E(u)|u|^{\beta (n-1)}e^{\sum_1^{\iy}t_i
%u^i}
%  e^{-\beta \sum_1^{\iy}\frac{u^{-i}}{i}
%\frac{\pl}{\pl t_i}}\tau_{n-1}(t)
\nonumber\\
&&\nonumber\\
 &=&\Big({\cal Y}_{}(t,s;
\rho)I(t,s)\Big)_n,
\eea
from which (2.3.9) follows.

Given the vertex operator $\BX_{12}(t,s;u,u^{-1})$, we now
 compute the corresponding Virasoro algebra, using (2.1.11),

\noindent$\displaystyle{
 u\frac{d}{du}u^k \BX_{12}(t,s;u,u^{-1})}$
\begin{eqnarray}
&=&\left( u^{k+1} \frac{d}{du}+ku^{k}\right)
 \BX_{12}(t,s;u,u^{-1})\nonumber\\
&=&\left. \left(u^{k+1}\frac{\pl}{\pl u}-v^{1-k}
 \frac{\pl}{\pl v}+ku^k\right) \BX_{12}
 (t,s;u,v) \right|_{v=u^{-1}}\nonumber\\
  &=&\left. \left(\frac{\pl}{\pl u}u^{k+1}-\frac{\pl}{\pl v}v^{1-k}
 -ku^k\right)\BX_{12}(t,s;u,v) \right|_{v=u^{-1}}\nonumber\\
 &=&\left. \left(\frac{\pl}{\pl u}u^{k+1}-\frac{\pl}{\pl v}v^{1-k}
 -k\theta u^k-k(1-\theta)v^{-k}\right)\BX_{12}(t,s;u,v)
   \right|_{v=u^{-1}}\nonumber\\
  &=&\Bigl[ {}^\beta\BJ_{k}^{(2)}(t)-{}^\beta\BJ_{-k}^{(2)}(-s)
-k\left( \theta ~{}^\beta\BJ_{k}^{(1)}(t)+(1-\theta)~~
 {}^\beta\BJ_{-k}^{(1)}(-s) \right),\nonumber\\
   & &\hspace{5cm}\BX_{12}(t,s;u,u^{-1})\Bigr]\Bigr|_{\beta=1}\nonumber\\
% \\ &=&\left[ \BJ_{k}^{(2)}(t)-\BJ_{-k}^{(2)}(-s)
%-k\left( \theta \BJ_{k}^{(1)}(t)+(1-\theta)\BJ_{k}^{(1)}(-s)
%   \right),\BX_{12}\right],\\
&=&\left[ {\cal V}_k^{(2)},\BX_{12}(t,s;u,u^{-1})\right],
   \end{eqnarray}
from which (2.3.6) follows. Verifying (2.3.7) goes by explicit
computation, using (5.0.2).

 Since (by virtue of (2.3.6)) for Lebesgue measure on $S^1$,
\bean \left[ \VR^{(2)}_k,  {\cal Y}(t,s,\rho =1)
 \right]&=&
\left[ \VR^{(2)}_k, \int_{S^1}  \BX_{12}(t,s;u,u^{-1})
 \frac{du}{2\pi i u} \right]\\&=&
  \int_{S^1} \left[\VR_k^{(2)}, \BX_{12}(t,s;u,u^{-1})
  \frac{du}{2\pi i u}  \right]\\&=&
   \int_{S^1}\frac{du}{2\pi i} ~\frac{d}{d u}u^{k+1}
   \frac{\BX_{12}(t,s;u,u^{-1})}{u} \\ & =&0,
\eean we have, using the notation (2.3.10) and
 the fact that, for $n\geq 0$ and $\rho=1$, the
integrals $I_n=n!\tau_n(t,s)$ are fixed points for
 ${\cal Y}(t,s,dz)$; hence
   \bean 0&=&\left([{\cal V}^{(2)}_{k},({\cal
Y}(t,s;dz))^n] I\right)_n\\ &=&\left({\cal
V}^{(2)}_{k} {\cal Y}(t,s;dz)^n I-{\cal
Y}(t,s;dz)^n{\cal V}^{(2)}_{k}I\right)_n\\
&=&\left({\cal V}^{(2)}_{k}I-{\cal Y}(t,s;dz)^n{\cal
V}^{(2)}_{k}I\right)_n. \eean Taking the $n^{{\rm
th}}$ component and taking into account the presence
of $\Lambda^{-1}$ in $\BX_{12}(t,s;u,u^{-1})$, we find
\bean 0&=&\left({\cal V}^{(2)}_{k}I-{\cal
Y}(t,s;dz)^n{\cal V}^{(2)}_{k}I\right)_{n}\\
 &=&{\cal
V}^{(2)}_{k}I_{n}-\int_{S_{1}}\frac{du}{2\pi
iu}e^{\sum_{1}^{\iy}(t_{i}u^i-s_{i}u^{-i})}
e^{-\sum_{1}^{\iy}\left(\frac{u^{-i}}{i} \frac{\pl}{\pl t_{i}}-
\frac{u^{i}}{i} \frac{\pl}{\pl s_{i}}\right)}
\\
&&   \hspace{2cm} \ldots \int_{S_{1}}\frac{du}{2\pi iu}
e^{\sum_{1}^{\iy}(t_{i}u^i-s_{i}u^{-i})}
e^{-\sum_{1}^{\iy}\left(\frac{u^{-i}}{i} \frac{\pl}{\pl t_{i}}-
\frac{u^{i}}{i} \frac{\pl}{\pl s_{i}}\right)} {\cal
V}^{(2)}_{k}I_{0}. \eean In the notation (2.3.5) and (5.0.4),
${\cal V}_{k}^{(2)}(t,s)$ has the following form \bean {\cal
V}^{(2)}_{k}(t,s)&=&\frac{1}{2}\left(J^{(2)}_{k}(t)-J^{(2)}_{-k}(-s)+
 (2n+k+1)J^{(1)}_{k}(t)-(2n-k+1)J^{(1)}_{-k}(-s)\right)\\ &
&\hspace{1cm}-k\left(\theta
J^{(1)}_{k}(t)+(1-\theta)J^{(1)}_{-k}(-s)\right) .
\eean

Working out the expression above leads to the expression written
out in (4.0.7) and one checks immediately that, given $\tau_0=1$,
$$ {\cal V}^{(2)}_{k}(t,s)\tau_0=0\quad\mbox{only for}~~
\left\{\begin{array}{l} k =-1 ,~~ \theta= 0\\ k=0 , ~~\theta~~
{arbitrary}\\
  k=1 ,~~
\theta= 1\end{array}\right\} , $$ ending the proof of
Theorem 2.5, except for the last equality of (2.3.4),
which follows easily from the last Proposition of
subsection 2.1. Indeed,
 $$
m_{\iy}(0,0)= \Bigl(\oint_{S^1}z^{k-\ell
}%e^{\sum_1^{\iy} (t_iz^i-s_iz^{-i})}
 \frac{dz}{2\pi
iz}\Bigr)_{1\leq k,\ell < \iy}=I_{\iy},
 $$
  where $I_{\iy}$
denotes the semi-infinite identity matrix. Finally,
one uses identity (2.1.27) and the fact that all the
determinants of submatrices $(I_{\iy})^{\lb\mu}$ (in
the notation (2.1.27)) are zero, except when the Young
diagrams $\lb$ and $\mu$ are equal. \qed

\remark  According to the strong Szeg\"o theorem, we
have: $$ \tau_n = \exp (-\sum_1^{\iy}
kt_ks_k)~~\mbox{for} ~~n \rightarrow \iy ,$$ provided
$\sum_1^{\iy} k(|t_k|^2+|s_k|^2)<\iy$. Therefore the
Toeplitz case yields boundary conditions for $\tau_n$
at both extremities, namely $n=0$ and $n=\iy$.
%\newpage

\subsection{Toeplitz matrices and the structure of $L_1$ and $L_2$}

The associated 2-Toda matrices $L_{1}$ and $L_{2}$ have a very
peculiar structure, when the initial $m_{\iy}$ matrix is
 Toeplitz, as we shall see in the main theorem
  of this section. Throughout, we shall be using the multiplication
  operator identity $z^{\top}=z^{-1}$ with regard to the
   inner
  product (2.3.1). This characterizes the Toeplitz case.
   Remember, from section 2.1, the
  polynomials (combining (2.1.4) and
  (2.1.24))
  $$
   p_n^{(1)}(z)=z^n\frac{\tau_n(t-[z^{-1}],s)}{\tau_n(t,s)}
   ~~\mbox{and}~~
   p_n^{(2)}(z)=z^n\frac{\tau_n(t,s+[z^{-1}])}{\tau_n(t,s)}
$$ are bi-orthogonal for the special inner-product
(2.3.1); also consider the vector notation:
 $$
p^{(i)}=(p_0^{(i)},p_1^{(i)},...)~,~~
 p_{\Lambda}^{(i)}=(p_1^{(i)},p_2^{(i)},...)
 ~~\mbox{and}~~h=\diag\left(\frac{\tau_1}{\tau_0},
  \frac{\tau_2}{\tau_1},...\right).
$$

\begin{theorem} The lower-triangular parts\footnote{In the formulae
 below $A_{-0}$ denotes the lower-triangular part of $A$, including the
diagonal.} of the matrices $L_{1}$ and $hL_2^{\top}h^{-1}$,
 arising in
 the context of a Toeplitz matrix $m_{\iy}$, are the projection of
 a rank 2 matrix:
\begin{eqnarray*}
L_{1}
 %&=&\sum^{-2}_{k=-\iy}\diag\left(\frac{p_{1-k}(\tilde\pl_t)
%\tau_{n+k+1}\circ\tau_{n}}{\tau_{n+k+1}\tau_{n}}\right)_{n\in\BZ}\Lambda^k+
%\left(\frac{\pl}{\pl t_1}\right)^2\log\tau_{n}\Lambda^{-1}+\frac{\pl}{\pl
%t_1}\log h_{n}\Lambda^0+\Lambda\\
&=&-\left(hp_{\Lambda}^{(1)}(0)\otimes h^{-1}p^{(2)}(0)\right)_{-0}+\Lambda\\
& & \\
  hL_2^{\top }h^{-1}
 %&=&\sum^{-2}_{k=-\iy}\diag\left(\frac{p_{1-k}(-\tilde\pl_s)
%\tau_{n+k+1}\circ\tau_{n}}{\tau_{n+k+1}\tau_{n}}\right)\Lambda^k+
%\left(\frac{\pl}{\pl s_1}\right)^2\log\tau_{n}\Lambda^{-1}-\frac{\pl}{\pl
%s_1}\log h_{n}\Lambda^0+\Lambda\\
&=&-\left(hp_{\Lambda}^{(2)}(0)\otimes h^{-1}p^{(1)}(0)\right)_{-0}+\Lambda.
\end{eqnarray*}

\end{theorem}

\vspace{0.6cm}

\begin{corollary} {\em (Unsymmetric identities)}
In particular\footnote{See footnote 8 for notation
$p_k(\tilde \pl_t)$ and $p_k(-\tilde \pl_s)$. The
bi-orthogonal polynomials $p_k^{(i)}$ should not be
confused with the Schur polynomials $p_k$.}
 $$
\begin{array}{ll}
p^{(1)}_{n+1}(0)p^{(2)}_{n+1}(0)=1-\frac{h_{n+1}}{h_{n}}&~~~p^{(2)}_{n+1}(0)p^{(1)}
_{n+1}(0)=1-
\frac{h_{n+1}}{h_n}\\
 \\
p^{(1)}_{n+1}(0)p^{(2)}_n(0)=-\frac{\pl}{\pl t_1}\log
h_{n}&~~~p^{(2)}_{n+1}(0)p^{(1)}_n(0)=\frac{\pl}{\pl s_{1}}\log h_n\\
 \\
p^{(1)}_{n+1}(0)p^{(2)}_{n-1}(0)=-\frac{h_{n-1}}{h_n}\left(\frac{\pl}{\pl
t_1}\right)^2\log
\tau_n&~~~p^{(2)}_{n+1}(0)p^{(1)}_{n-1}(0)=-\frac{h_{n-1}}{h_n}\left(\frac{\pl}{\pl
s_1}\right)^2\log \tau_n
\end{array}
$$
\bea
p^{(1)}_{n+1}(0)p^{(2)}_{n-k}(0)&=&
 -\frac{h_{n-k}}{h_n}\frac{p_{k+1}(\tilde\pl_t)
\tau_{n-k+1}\circ\tau_n}
{\tau_{n-k+1}\tau_n}  \nonumber\\
p^{(2)}_{n+1}(0)p^{(1)}_{n-k}(0)&=&
 -\frac{h_{n-k}}{h_n}\frac{p_{k+1}(-\tilde\pl_s)
\tau_{n-k+1}\circ\tau_n}
{\tau_{n-k+1}\tau_n} ,~~k \geq 0.
\eea
\end{corollary}

\begin{corollary} {\em (Symmetrized identities)}
We also have for $n>m$,
 \bean &&
\left(\frac{h_{n}}{h_{m+1}}\right)^{2}\left(1-\frac{h_{n+1}}{h_{n}}\right)
\left(1-\frac{h_{m+1}}{h_{m}}\right)\\ &&  \\&&
~~~~=\frac{1}{\tau_{m+2}^2\tau_n^2}
 \left(p_{n-m}(\tilde\pl_{t})\tau_{m+2}\circ
\tau_{n}\right)~.~\left(p_{n-m}(-\tilde\pl_{s})\tau_{m+2}
\circ\tau_{n}\right)
~.
\eean

 \medbreak
\noindent In particular, for $m=n-1$,
\be
\left(1-\frac{h_{n+1}}{h_{n}}\right)\left(1-\frac{h_{n}}
{h_{n-1}}\right)
=-\frac{\pl}{\pl t_{1}}\log h_{n}\frac{\pl}{\pl s_{1}}\log h_{n}.
\ee

\end{corollary}

Identity (2.4.2) was already observed by Hisakado in \cite{H}. We
first need a Lemma, which explains the peculiar structure of the
bi-orthogonal polynomials $p_n^{(1)}(y)$ and $p_n^{(2)}(z)$,
associated with the inner-product (2.3.1):

\begin{lemma} (Hisakado \cite{H}) \, The following holds:
\bea
p^{(1)}_{n+1}(z)-zp_n^{(1)}(z)&=&p^{(1)}_{n+1}(0)z^np_n^{(2)}(z^{-1})\nonumber\\
p^{(2)}_{n+1}(z)-zp_n^{(2)}(z)&=&p^{(2)}_{n+1}(0)z^np_n^{(1)}(z^{-1}).
\eea
\end{lemma}

\proof The following orthogonality relations hold for $1\leq i\leq n$:
$$
\la p_{n+1}^{(1)}(z)-zp_n^{(1)}(z),z^i\ra = \la
p_{n+1}^{(1)}(z),z^i\ra -\la p_n^{(1)}(z),z^{i-1}\ra=0
$$
and
$$
\la z^np^{(2)}_{n}(z^{-1}),z^i\ra =\la z^{n-i},p_{n}^{(2)}(z)\ra =0.
$$
Therefore the two $n^{th}$ degree polynomials
$p_{n+1}^{(1)}(z)-zp_n^{(1)}(z)$ and  $z^np^{(2)}_{n}(z^{-1})$ must
be proportional and since
$$
p^{(1)}_{n+1}(z)-zp_n^{(1)}(z)\Big|_{z=0}=p^{(1)}_{n+1}(0),~~\mbox{and}~~z^np^{(
2)}_{n}
(z^{-1})\Big|_{z=0}=1,
$$
the first identity (2.4.3) follows. The second one follows by
duality.\qed

\noindent\underline{\sl Proof of Theorem 2.7 and Corollaries 2.8 and 2.9}:
On the one
hand\footnote{Define $p_{-1}^{(2)}(z)=0$.},

 $\la
p_{n+1}^{(1)}(z)-zp_{n}^{(1)}(z),p^{(2)}_{m+1}(z)
 -zp^{(2)}_{m}(z)\ra\hfill
n>m\geq -1$ \bea &=&-\la
zp_{n}^{(1)}(z),p_{m+1}^{(2)}(z)\ra\nonumber\\ &=&-\la
p_{n+1}^{(1)}(z)+\ldots
+(L_1)_{n,m+1}p_{m+1}^{(1)}(z)+\ldots~~,p_{m+1}^{(2)}(z)\ra\nonumber\\
&=&-(L_1)_{n,m+1}\la
p^{(1)}_{m+1}(z),p^{(2)}_{m+1}(z)\ra\nonumber\\
&=&-(L_1)_{n,m+1}h_{m+1}, \eea and, on the other hand,

 $\la
p_{n+1}^{(1)}(z)-zp_{n}^{(1)}(z),p^{(2)}_{m+1}(z)-zp^{(2)}_{m}(z)\ra\hfill
n\geq m\geq -1$
\bea
&=&\la p_{n+1}^{(1)}(0)z^np_n^{(2)}(z^{-1}),p_{m+1}^{(2)}(0)
z^mp_m^{(1)}(z^{-1})\ra\nonumber\\
&=&p_{n+1}^{(1)}(0)p^{(2)}_{m+1}(0)\la
z^{n-m}p^{(1)}_m(z),p^{(2)}_n(z)\ra\nonumber\\ &=&
p^{(1)}_{n+1}(0)p^{(2)}_{m+1}(0)\la
p^{(1)}_{n}(z)+\ldots,p_{n}^{(2)}(z)\ra\nonumber\\ &=&
p^{(1)}_{n+1}(0)p^{(2)}_{m+1}(0)h_{n}.
\eea
Comparing (2.4.4) and (2.4.5) yields
\be
(L_{1})_{n,m+1}=-h_{n}p_{n+1}^{(1)}(0)h^{-1}_{m+1}
 p_{m+1}^{(2)}(0), ~~n>m\geq -1,
\ee proving the first expression of Theorem 2.7; the
second one is obtained by the usual duality
$L_1\mapsto hL_2^{\top}h^{-1},~ t\leftrightarrow -s$
 and so
$p^{(1)} \leftrightarrow p^{(2)}$ (see formulae in the
beginning of this section).
 For
$n=m$, we compute

$\la
p_{n+1}^{(1)}(z)-zp_{n}^{(1)}(z),p^{(2)}_{n+1}(z)-zp^{(2)}_{n}(z)\ra$
\bean &=&\la p^{(1)}_{n+1}(z),p_{n+1}^{(2)}(z)\ra +\la
zp^{(1)}_{n}(z),zp_n^{(2)}(z)\ra -\la
zp_{n}^{(1)}(z),p_{n+1}^{(2)}(z) \ra\\ & &\hspace{4cm} -\,\la
p^{(1)}_{n+1}(z),zp^{(2)}_{n}(z)\ra\\
&=&h_{n+1}+h_{n}-h_{n+1}-h_{n+1}\\ &=&h_{n}-h_{n+1}, \eean which
upon comparison with (2.4.5) for $n=m$ yields the first line of
corollary 2.8:
\be
p^{(1)}_{n+1}(0)p^{(2)}_{n+1}(0)h_{n}=h_{n}-h_{n+1}. \ee Remember
from (2.1.9) and (2.1.12), we have:
\begin{eqnarray*}
L_{1}
 &=&\sum^{-2}_{k=-\iy}\diag\left(\frac{p_{1-k}(\tilde\pl_t)
\tau_{n+k+1}\circ\tau_{n}}{\tau_{n+k+1}\tau_{n}}
 \right)_{n\in\BZ}\Lambda^k\\ && \hspace{2cm}+
\left(\frac{\pl}{\pl t_1}\right)^2\log\tau_{n}~\Lambda^{-1}+\frac{\pl}{\pl
t_1}\log h_{n}~\Lambda^0+\Lambda\\
%&=&-\left(hp_{\Lambda}^{(1)}(0)\otimes
%h^{-1}p^{(2)}(0)\right)_{-0}+\Lambda\\
 & & \\
  hL_2^{\top }h^{-1}
 &=&\sum^{-2}_{k=-\iy}\diag\left(\frac{p_{1-k}(-\tilde\pl_s)
\tau_{n+k+1}\circ\tau_{n}}{\tau_{n+k+1}\tau_{n}}\right)
 \Lambda^k  \\ && \hspace{2cm} +
\left(\frac{\pl}{\pl s_1}\right)^2\log\tau_{n}~\Lambda^{-1}-\frac{\pl}{\pl
s_1}\log h_{n}~\Lambda^0+\Lambda. \\
%&=&-\left(hp_{\Lambda}^{(2)}(0)\otimes
%h^{-1}p^{(1)}(0)\right)_{-0}+\Lambda;
\end{eqnarray*}
Together with the theorem, this yields corollary 2.8.

Finally, upon multiplying relations (2.4.1), setting
 $m+1=n-k$, $n>m$, and using the relation above, one obtains,
 using Corollary 2.8,
\bea \lefteqn{\left(\frac{h_{m+1}}{h_{n}}\right)^2
\frac{\left(p_{n-m}(\tilde\pl_{t})\tau_{m+2} \circ
\tau_{n}\right)\,\,\left(p_{n-m}
(-\tilde\pl_{s})\tau_{m+2}\circ
\tau_{n}\right)}{\tau^2_{m+2}\tau^2_{n}}}\nonumber\\
&&\nonumber\\
 &=&p^{(1)}_{n+1}(0)p^{(2)}_{m+1}(0)p^{(2)}_{n+1}(0)p^{(1)}_{m+1}(0)\nonumber\\
&=&\left(1-\frac{h_{n+1}}{h_{n}}\right)
\left(1-\frac{h_{m+1}}{h_{m}}\right)\nonumber\\
&=&\frac{1}{\tau^2_{n+1}\tau^2_{m+1}}(\tau^{2}_{n+1}-\tau_{n}\tau_{n+2})(\tau^{2
}_{m+1}-\tau_{m }\tau_{m+2}), \eea which is precisely
Corollary 2.9. Relation (2.4.2) is a special case of
(2.4.8), setting $m=n-1$. \qed

\medskip\noindent{\it Proof of Theorem 0.3:\/} The
structure of $L_1$ and $L_2$ follows from theorem 2.7.
The statement about the mathematical expectation
follows from:
%\bigbreak
\bean
  \lefteqn{p_n^{(1)} (t,s;z)}\\
  &=&z^n \frac{\tau_n(t-[z^{-1}],s)}{\tau_n(t,s)}\\
  &=& \sum_{k=0}^{n}z^k~\frac{p_{n-k}(-\tilde \pl _t)\tau_n(t,s)}
  {\tau_n(t,s)}\\
  &=&\frac{1}{\tau_n} \sum_{k=0}^{n}z^k~
  \int_{U(n)}p_{n-k}(-\Tr M,-\frac{1}{2} \Tr M^2,
  -\frac{1}{3} \Tr M^3,...)e^{\sum_1^{\iy}\Tr (t_iM^i-s_i\bar
M^i)} dM,
%\\
%    &=&  \sum_{k=0}^{n}z^k~ (-1)^{n-k}E_{U(n)}
%     s_{\lambda_{n-k}}( M), ~~\mbox{with }~ \lambda_i=
%     \overbrace {1,...,1}^i,
\eean
 and similarly
 \bean
 \lefteqn{p_n^{(2)} (t,s;z)}\\
  &=&z^n \frac{\tau_n(t,s+[z^{-1}])}{\tau_n(t,s)}\\
  &=& \sum_{k=0}^{n}z^k~\frac{p_{n-k}(\tilde \pl _s)\tau_n(t,s)}
  {\tau_n(t,s)}\\
  &=&\frac{1}{\tau_n} \sum_{k=0}^{n}z^k~
  \int_{U(n)}p_{n-k}(-\Tr \bar M,-\frac{1}{2} \Tr \bar M^2,
  -\frac{1}{3} \Tr \bar M^3,...)e^{\sum_1^{\iy}\Tr (t_iM^i-s_i\bar
M^i)} dM.
%\\
%    &=&   \sum_{k=0}^{n}z^k~(-1)^{n-k} E_{U(n)}
%     s_{\lambda_{n-k}}(\bar M), ~~\mbox{with }~ \lambda_i=
%     \overbrace{1,...,1}^i.
\eean
 Finally, we check the Hamiltonian flow
statement for the first flow. Indeed, from the
equations for $\Psi$ (after (2.1.7)), (2.1.23),
Theorem 2.7, and the first relation of Corollary 2.8,
 it follows that ($h_{-1}=0$)
 \bean
\frac{\pl x_n}{\pl t_1}
 &=& \left. \frac{\pl p_n^{(1)}(t,s;z)}{\pl t_1}\right|_{z=0}\\
 &=& \left. -\left((L_1)_-p^{(1)}\right)_n \right|_{z=0}\\
 &=& h_n p^{(1)}_{n+1}(t,s;0)\sum_{i=0}^{n-1} \frac
 { p^{(1)}_{i}(t,s;0) p^{(2)}_{i}(t,s;0)}{h_i}\\
 &=& h_n x_{n+1} \sum_{i=0}^{n-1} \frac
 { x_{i} y_{i}}{h_i}\\
 &=&  h_n x_{n+1} \sum_{i=0}^{n-1}\left(
 \frac{1}{h_i}-\frac{1}{h_{i-1}}\right)
% \\ &=&
 = x_{n+1}\frac{h_n}{h_{n-1}}
% \\ &=&
  =x_{n+1}  (1-x_ny_n),
 \eean
 and similarly for the other coordinates
ending the proof of Theorem 0.3.\qed

\section{Painlev\'e equations for $O(n)$ and $Sp(n)$ integrals}

\subsection{Painlev\'e equations associated with the Jacobi weight}
%Remember  the definition (1.0.1)
% of the Jacobi weight, namely $\rho_{\alpha\beta}(z)dz:=(1-z)^{\alpha}(1+z)^{\beta} dz$;
%We now state:

\begin{theorem}{\em (Painlev\'e equation and the Jacobi weight)}
~The function $H_n(x)=x\frac{d}{dx}\log \tau_n(x)$,
with
\be
\tau_n(x):=c_n
 \int_{[-1,1]^n}\Delta_n(z)^2
 \prod_{k=1}^{n} e^{x z_k}
  (1-z_k)^{\alpha}(1+z_k)^{\beta} dz_k,
 \ee
satisfies the Painlev\'e V equation $(a:=\alpha+\beta,
~b:=\alpha -\beta,$ and $\alpha,\beta > -1)$
 \bea
   && \hspace{-.4cm}
x^2 H^{\prime\prime\prime}+xH^{\prime\prime}
 +6x{H^{\prime} }^2-\left(4H+4x^2-4bx
  +(2n+a)^2\right)H^{\prime}+(4x-2b)H\nonumber\\
&&\hspace{7cm}+2n(n+a)x-bn(2n+a)=0,\nonumber\\ \eea
%for $a:=\alpha+\beta, ~b:=\alpha -\beta,$ and
%$\alpha,\beta > -1$.
with initial condition
\be
H(0)=0~~\mbox{and}~~H^{\prime}(0)%=\left.\frac{d}{dx}
%\log \tau_n(x)\right|_{x=0}
=\frac{-nb} {a+2n}. \ee
\end{theorem}

\begin{corollary}

\be
\tilde H_n(x)=x\frac{d}{dx}\log e^{-cx}\tau_{n}(2x)
\ee satisfies the Painlev\'e V equation \bea
&&\hspace{-0.7cm}\frac{1}{2}x^2{\tilde
H}^{\prime\prime\prime}
 +\frac{1}{2}x{\tilde H}^{\prime\prime}
 +3x({{\tilde H}^{\prime} })^2
  -\frac{1}{2}\left(4\tilde H+16x^2-8(b+c)x
  +(2n+a)^2\right){\tilde H}^{\prime}+\nonumber\\
&&\hspace{-0.7cm}(8x-2(b+c))\tilde
H+(4n(n+a)+c(2b+c))x-\frac{1}{2}(2n+a)(2n(b+c)+ac)=0.\nonumber\\
\eea with
\be
\tilde H(0)=0~~\mbox{and}~~\tilde H^{\prime}(0)=
%\left.\frac{d}{dx}\log e^{-cx} \tau_n(2x)\right|_{x=0}=
 -\frac{2n(b+c)+ac} {2n+a}. \ee

\end{corollary}

\noindent\underline{\sl Proof of Theorem 3.1}:
Since $\alpha,\beta
>-1$, the boundary condition (0.0.6) on $\rho(z)$ is fulfilled; so,
we may apply Theorem 2.4. We set $$
a_0=1,a_1=0,a_2=-1,b_0=\alpha
-\beta=:b,b_1=\alpha+\beta=:a $$ and all other
$a_i=b_j=0$ in (2.2.12), implying that
 $$
 \tau_n(t_1,t_2,...):=
 \int_{[-1,1]^n}\Delta_n(z)^2
 \prod_{k=1}^{n} e^{\sum_1^{\iy} t_i z_k^i }
  (1-z_k)^{\alpha}(1+z_k)^{\beta} dz_k
$$ satisfies the equations %\linebreak
($m=1,2,3,\ldots$)\footnote{The $J_{m}^{(i)}$ below
are the ones of (5.0.3), for $\beta =2$.}:
 \bean 0
&=&\JR_{m-2}^{(2)}\tau_{n}\\ &=&\sum_{k\geq 0}\left.
\left(-a_{k}\sum_{i+j=k+m-2}:\,^{\beta}\BJ^{(1)}_{i}
 \,^{\beta}\BJ^{(1)}_{j}:+b_{k}
\,^{\beta}\BJ^{(1)}_{k+m-1}\right)\right|_{\beta=2} \tau_{n}\\
&=&\left(J_m^{(2)}-J_{m-2}^{(2)}-2n J_{m-2}^{(1)}
 +(2n+a) J_m^{(1)}+b
J_{m-1}^{(1)}-n^2\dt_{m,2}+nb\dt_{m,1}\right)
\tau_n.
\eean
Then introducing the function $F_n:=\log \tau_n(t)$, the two first
Virasoro constraints for $m=1,2$ divided by $\tau_n$ are
given by
\bea
\frac{{\cal J}^{(2)}_{-1}\tau_n}{\tau_n}&=&\left(\sum_{i\geq 1}it_i\frac{\pl}{\pl
t_{i+1}}-
\sum_{i\geq 2}it_i\frac{\pl}{\pl t_{i-1}}
+(2n+a)\frac{\pl}{\pl t_1}\right) F_n+n(b-t_1)=0
\nonumber\\
 \frac{{\cal J}^{(2)}_{0}\tau_n}{\tau_n}&=& \left(\sum_{i\geq
1}it_i\frac{\pl}{\pl
t_{i+2}}-
\sum_{i\geq 1}it_i\frac{\pl}{\pl t_{i}} +
 b\frac{\pl}{\pl t_1}
 +\frac{\pl^2}{\pl t_1^2}
+(2n+a)\frac{\pl}{\pl t_2}\right) F_n \nonumber\\
 &&\hspace{7cm}+\left(\frac{\pl F_n}{\pl t_1}\right)^2-n^2=0.\nonumber\\
\eea These expressions and their first $t_1$- and $t_2$-
derivatives, evaluated along the locus $$ \LR:=\{t_1=x, ~\mbox{all
other}~t_i=0 \} $$ read as follows:
 \bean
 0&=&\left.\frac{{\cal
J}^{(2)}_{-1}\tau_{n}}{\tau_{n}}\right|_{\LR}=\left.\left(t_{1}\frac{\pl}{\pl
t_{2}}+ (2n+a)\frac{\pl}{\pl
t_{1}}\right)F_n+n(b-t_{1})\right|_{\LR}\\
 0&=&\left.\frac{{\cal
J}^{(2)}_{0}\tau_{n}}{\tau_{n}}\right|_{\LR}=
\left(t_{1}\frac{\pl}{\pl t_{3}}+(b-t_1)
 \frac{\pl}{\pl
t_{1}}+(2n+a)\frac{\pl}{\pl t_2}+\frac{\pl^2}{\pl
t^2_1}\right)F_n\\ & &\hspace{5cm}+\left.\left(\frac{\pl F_n}{\pl
t_1}\right)^2-n^2\right|_{\LR}\\ & & \\ 0&=&
 \left.\frac{\pl}{\pl t_1} \frac{{\cal
J}^{(2)}_{-1}\tau_{n}}{\tau_{n}}\right|_{\LR}
 =\left(\sum_{i\geq
 1}it_{i}\frac{\pl^{2}}{\pl t_{i+1}\pl t_{1}}
 +\frac{\pl}{\pl t_{2}}-
\sum_{i\geq 2}it_{i} \frac{\pl^{2}}{\pl t_{i-1}\pl t_{1}}\right.\\
& &\left.\left.\hspace{4cm}+(2n+a)\frac{\pl^2}{\pl
t_1^2}\right)F_n\right|_{\LR}-n\\ && \hspace{2.3cm}
=\left.\left(t_1\frac{\pl^2}{\pl t_2\pl t_1}+\frac{\pl}{\pl
t_{2}}+(2n+a)\frac{\pl^2}{\pl t_1^2}\right)F_n
 \right|_{\LR}-n\\ & &
\\
 0&=&\left.\frac{\pl}{\pl t_1}\frac{{\cal
J}^{(2)}_{0}\tau_{n}}{\tau_{n}}\right|_{\LR}=\left(\sum_{i\geq
1}it_{i}\frac{\pl^{2}}{\pl t_{i+2}\pl
t_{1}}+\frac{\pl}{\pl t_{3}}- \sum_{i\geq
1}it_{i}\frac{\pl^{2}}{\pl t_{i}\pl
t_{1}}-\frac{\pl}{\pl t_1}+b\frac{\pl^2}{\pl
t_1^2}\right.\\ & &\hspace{3cm}\left.\left.
+\frac{\pl^3}{\pl t_1^3} +(2n+a)\frac{\pl^2}{\pl
t_2\pl t_1}\right)F_n+2\frac{\pl F_n}{\pl
t_{1}}\frac{\pl^2 F_n}{\pl t_1^2}\right|_{\LR}\\
&&\hspace{2.3cm}=\left(t_1\frac{\pl^2}{\pl t_3\pl
t_1}+\frac{\pl}{\pl t_3}+(b-t_1)\frac{\pl^2}{\pl
t_1^2} -\frac{\pl}{\pl t_1}+\frac{\pl^3}{\pl
t_1^3}\right.\\ &
&\hspace{3cm}\left.+(2n+a)\frac{\pl^2}{\pl t_2\pl
t_1}\right)F_n+2\left.\frac{\pl F_n}{\pl
t_1}\frac{\pl^2 F_n}{\pl t_1^2}\right|_{\LR}\\ & & \\
0&=&\left.\frac{\pl}{\pl t_2}\frac{{\cal
J}^{(2)}_{-1}\tau_{n}}{\tau_{n}}\right|_{\LR}=
\left(\sum_{i\geq 1}it_{i}\frac{\pl^{2}}{\pl
t_{i+1}\pl t_{2}}+2\frac{\pl}{\pl t_{3}}- \sum_{i\geq
2}it_{i}\frac{\pl^{2}}{\pl t_{i-1}\pl
t_{2}}-2\frac{\pl}{\pl t_1}\right.\\ &
&\hspace{4.5cm}+(2n+a)\left.\left.\frac{\pl^2}{\pl
t_1\pl t_2}\right)F_n\right|_{\LR}\\
&&\hspace{2.4cm}=\left(t_1\frac{\pl^2}{\pl
t_2^2}+2\frac{\pl}{\pl t_3}-2\frac{\pl}{\pl t_1}+
(2n+a)\frac{\pl^2}{\pl t_1\pl t_2}\right)F_n. \eean
The five equations above form a (triangular) linear
system in five unknowns $$ \left.\frac{\pl F_n}{\pl
t_2}\right|_{\LR},
 \quad\left.\frac{\pl F_n}{\pl t_3}\right|_{\LR},\quad\left.
 \frac{\pl^2F_n}{\pl
t_1\pl t_2}\right|_{\LR},\quad\left.\frac{\pl^2F_n}{\pl t_1\pl
t_3}\right|_{\LR},\quad \left.\frac{\pl^2F_n}{\pl
t_2^2}\right|_{\LR}.\quad $$
 Setting $t_1=x$ and $F'_n=\pl F_n/ \pl x$, the solution is
given by the following expressions, \bean \left.\frac{\pl F_n}{\pl
t_2}\right|_{\LR}&=&-\frac{1}{x}
\Bigl((2n+a)F^{\prime}_n+n(b-x)\Bigr)\\ & & \\ \left.\frac{\pl
F_n}{\pl
t_3}\right|_{\LR}&=&-\frac{1}{x^2}\Bigl(x\left(F^{\prime\prime}_n+
F_n^{\prime 2}+(b-x)F^{\prime}_n+n(n+a)\right)-(2n+a)
\left((2n+a)F^{\prime}_n+bn\right)\Bigr)\\ & & \\
\left.\frac{\pl^2 F_n}{\pl t_1\pl
t_2}\right|_{\LR}&=&-\frac{1}{x^2}
 \Bigl((2n+a)(xF^{\prime\prime}_n-F^{\prime}_n)-bn\Bigr)
\\
\left.\frac{\pl^2 F_n}{\pl t_1\pl t_3}\right|_{\LR}
 &=&-\frac{1}{x^3}\Bigl(
x^2(F^{\prime\prime\prime}_n+2
F_n^{\prime}F_n^{\prime\prime})-x\left((x^2-bx+1)
F_n^{\prime\prime}+F_n^{\prime
2}+bF_n^{\prime}+(2n+a)^2F_n^{\prime\prime}\right.\\
 &&\hspace{3cm} \left.
+n(n+a)\right)+2(2n+a)^2F_n^{\prime}+2bn(2n+a)\Bigr)\\ & & \\
\left. \frac{\pl^2 F_n}{\pl t_2^2}\right|_{\LR}
 &=&\frac{1}{x^3}\Bigl(
x\left(2F^{\prime 2}_n+2b
F^{\prime}_n+((2n+a)^2+2)F^{\prime\prime}_n+2n(n+a)\right)\\
 & &\hspace{4cm} -3(2n+a)^2F^{\prime}_n
  -3bn(2n+a)\Bigr).
\eean

Putting these expressions into the KP-equation (Theorem 2.1), and setting
$$
G(x):=F_n'(x)=\frac{d}{dx} \log \tau_n(x),
$$
we find
\bea
&&x^3G^{\prime\prime\prime}+4x^2G^{\prime\prime}
+x\left(-4x^2+4bx+2-(2n+a)^2\right)G'
+8x^2GG^{\prime}+6x^3{G^{\prime}}^2\nonumber\\
&&+2xG^2+\left(2bx-(2n+a)^2\right)G+n(2x-b)(n+a)-bn^2=0.
\nonumber\\
 \eea
Finally, the function
$$H(x):=xG(x)=x\frac{d}{dx} \log \tau_n(x)$$ satisfies
\bea
&&x^2H^{\prime\prime\prime}+xH^{\prime\prime}
 +6x{H^{\prime} }^2-\left(4H+4x^2-4bx
  +(2n+a)^2\right)H^{\prime}+(4x-2b)H\nonumber\\
&&\hspace{7cm}+2n(n+a)x-bn(2n+a)=0.\nonumber\\ \eea According to
Cosgrove \cite{C}, this 3rd order equation can be transformed into
a master Painlev\'e equation, which one recognizes to be
Painlev\'e V; see Appendix 2.

From section 8, Appendix 4, identity (8.0.5), it now
follows that $$H_n'(0)=
\frac{\tau'_{n}(0)}{\tau_{n}(0)}=
 \displaystyle{\frac{\sum^n_{i=1}\int_{[-1,1]^n}
\Delta_{n}
(z)^{2}z_{i}\prod^n_{k=1}\rho_{(\alpha,\beta)}(z_k)dz_k}{\int_{[-1,1]^n}
\Delta_{n}
(z)^{2}\prod^n_{k=1}\rho_{(\alpha,\beta)}(z_k)dz_k}}
%=n \frac{\beta -\alpha}{\beta+\alpha +2n}
=n \la y_1 \ra=\frac{-nb}{a+2n}. $$ This
ends the proof of Theorem 3.1.\qed

\medbreak

\noindent\underline{\sl Proof of Corollary 3.2}:  The differential
equation for $$ \tilde H(x)=x\frac{d}{dx}\log
e^{-cx}\tau_{n}(2x)=H(2x)-cx $$ is obtained by first setting
$x\mapsto 2x$ in the differential equation (3.1.2) and then setting
$H(2x)=\tilde H(x)+cx$. This leads to the differential equation
(3.1.5), which is, of course, also Painlev\'e V. Relation (3.1.6)
follows at once from (3.1.4). \qed

\subsection{Proof of Theorem 0.1 ($O(n)$ and $Sp(n)$)}

We give here a more detailed version of Theorem 0.1
(i):

\begin{proposition} Given the integral\footnote{ The integral over
the symplectic group $Sp(n-1)$ can be identified with
$O(2n)_-$. }
  $$
  I^{\pm}_{\ell}(x)=\int_{O_{\pm}(\ell)}
e^{x \Tr M }dM , $$ the expressions\footnote{In this
statement,
 we use the following notation:
$$[n]_{\mbox{\tiny{even}}}:=\max~
\{\mbox{even $x$, such that $x\leq n  \}$.}$$}
 $$ q_{\ell}(x)=\log e_{\ell}^{\pm}~
\frac{I^{\pm}_{\ell+2}}{I^{\pm}_{\ell}},~~\mbox{ with}~~ e^{+}_{\ell}=
 \frac{2}{[\ell+2]_{\mbox{\tiny{even}}}}~~
 \mbox{ and }~e^{-}_{\ell}=
 \frac{2}{[\ell+1]_{\mbox{\tiny{even}}}}
 %  ~e^{\pm}_{\ell}=
 %\frac{2}{[\ell+{2 \atop 1}]_{\mbox{\tiny{even}}}}
 $$
 satisfy the standard Toda lattice equations:
  $$ \frac{1}{4}\frac{\pl^2 q_{\ell} }{\pl x^2}  =-e^{q_{\ell}-q_{\ell
-1}}+ e^{q_{\ell +1}-q_{\ell }}.
$$
\end{proposition}

\begin{proposition} The function
\be
f^{\pm}_{\ell}(x) = x\frac{d}{dx}\log
\int_{O(\ell+1)_{\pm}~\mbox{\tiny{or}}
 ~Sp(\frac{\ell-1}{2})}e^{x\Tr M}dM
\ee is the \underline{unique} solution to the 3rd
order equation (i) in Theorem 0.1:
 $$
\mbox{\bf} \left\{\begin{array}{l}
\displaystyle{f^{\prime\prime\prime}+
 \frac{1}{x}
f^{\prime\prime}
 +\frac{6}{x}{f^{\prime} }^2-\frac{4}{x^2} f f^{\prime}
 -\frac{16 x^2+\ell^2}{x^2}f^{\prime}
+\frac{16}{x}f  +\frac{2(\ell^2-1)}{x}=0} \\  \\
\displaystyle{\mbox{with}~~  f_{\ell}^{\pm}(x)=x^2 \pm
\frac{x^{\ell+1}}{\ell !}+O(x^{\ell+2}),~\mbox{near}~
x=0 .}
\end{array}
 \right. $$ \vspace{-1.4cm}\be \ee
  This 3rd order equation can be transformed into
  the following second order equation in
$f$, quadratic in $f''$: \bean \frac{x^2}{4}
f^{\prime\prime 2}&=&-\left(xf^{\prime
2}-\left(4x^{2}+\frac{\ell^2}{4}\right)f^{\prime}+x(\ell^2-1)\right)f^{\prime}\\
& &\hspace{2cm}+(f^{\prime 2}-8xf^{\prime}+\ell^2-1)f+
4f^2, \eean which in turn leads to the standard
Painlev\'e equation (6.0.3) for $$ \alpha=-\beta =
\frac{(\ell+1)^2}{8},~ \gamma=0,~ \dt=-8 .$$

\end{proposition}

\medskip\noindent{\it Proof of Proposition 3.3:\/}
Proposition 1.1 and identity (2.2.11) imply

\begin{itemize}
  \item $I^+_{2n+1}(x)=n! e^x \tau_n(2x,0,...)$
  \item $I^+_{2n}(x)=n!  \tau_n(2x,0,...)$
  \item $I^-_{2n+1}(x)=n! e^{-x} \tau_n(2x,0,...)$
  \item $I^-_{2n}(x)=(n-1)!  \tau_{n-1}(2x,0,...)$.
\end{itemize}
Note that, since the functions $\tau_n(t-s)=\tau_n(t,s)$
 satisfy differential equation (ii) of Theorem 2.1,
   we obtain for the function $\tau_n(t)$, by subtracting
 two consecutive equations:
$$
\frac{\pl^2}{\pl
t_1^2}\log\frac{\tau_{n+1}}{\tau_n}= \frac{\tau_{n}\tau_{n+2}}{\tau_{n+1}^2}
-\frac{\tau_{n-1}\tau_{n+1}}{\tau_n^2},
$$
from which the standard Toda lattice equations follow. \qed

\medskip\noindent{\it Proof of Proposition 3.4:\/} From Corollary 3.2 it follows that

\be
\tilde H_n(x)=x\frac{d}{dx} \log \left( e^{-cx}
 \int_{[-1,1]^n}\Delta_n(z)^2
 \prod_{k=1}^{n} e^{2x z_k}
  %\rho_{(\alpha,\beta)}(z_k)
  (1-z_k)^{\alpha}(1+z_k)^{\beta}dz_k \right)
 \ee
satisfies the Painlev\'e V equation (3.1.5). Then in view of
Theorem 1.1, $\tilde H_n(x)$ corresponds to $f_{\ell}(x)$
in (3.2.1), when
the parameters $n$, $a=\alpha +\beta,~b=\alpha-\beta$
and $c$ take on the following values:

\vspace{0.5cm}

\noindent $O(\ell+1)_-,~\mbox{with}~\ell$ even:
 $n=\ell/2,~a=0,~ b=-1,~c=1$

\noindent $O(\ell+1)_-,~\mbox{with}~\ell$ odd:
 $n=(\ell-1)/2,~a=1,~ b=0,~c=0$

\noindent $O(\ell+1)_+,~\mbox{with}~\ell$ even:
 $n=\ell/2,~a=0,~ b=1,~c=-1$

\noindent $O(\ell+1)_+,~\mbox{with}~\ell$ odd:
 $n=(\ell +1)/2,~a=-1,~ b=0,~c=0$

\noindent $Sp(\frac{\ell-1}{2}),~\mbox{with}~\ell$ odd:
 $n=(\ell -1)/2,~a=1,~ b=0,~c=0$.

\vspace{0.4cm}

Setting these values into equation (3.1.5) leads at
once to equation (i) of Theorem 0.1, namely $$
f^{\prime\prime\prime}+ \frac{1}{x} f^{\prime\prime}
 +\frac{6}{x}{f^{\prime} }^2-\frac{4}{x^2} f f^{\prime}
 -\frac{16 x^2+\ell^2}{x^2}f^{\prime}
+\frac{16}{x}f +\frac{2(\ell^2-1)}{x}=0.
 $$
 Moreover,
for these values , we have that $b+c=ac=0$ and so from
(3.1.6), it follows that
 \be f_{\ell}(0)=\tilde H_{\ell}(0)=0 ~~\mbox{and}
  ~~f_{\ell}'(0)=\tilde H_{\ell}'(0)
 =- \frac{2n(b+c)+ac} {2n+a}=0
 .
\ee

According to appendix 2 (see Cosgrove \cite{C}), this
third order equation has a first integral, which is
second order in $f$ and quadratic in
$f^{\prime\prime}$, thus introducing a constant $c$:
\bean \frac{x^2}{4} f^{\prime\prime
2}&=&-\left(xf^{\prime
2}-\left(4x^{2}+\frac{\ell^2}{4}\right)f^{\prime}+x(\ell^2-1)\right)f^{\prime}\\
& &\hspace{2cm}+(f^{\prime
2}-8xf^{\prime}+\ell^2-1)f+4f^2-\frac{c}{4}. \eean
Evaluating this differential equation at $x=0$ leads
to, since $f(0)=0$,
 $$
c=\ell^2f^{\prime}(0)^2=0,\quad\mbox{using (3.2.4).}
 $$
  Setting $f=\bar f-\ell^2/4$ in order to
get the equation in Cosgrove's form \cite{CS},
  \bean
\frac{x^2}{4} \bar f^{\prime\prime 2}&=&-\left(x\bar
f^{\prime 2}-4x^{2}\bar
f^{\prime}-x(\ell^2+1)\right)\bar f^{\prime}\\ &
&\hspace{2cm}+(\bar f^{\prime 2}-8x\bar
f^{\prime}-(\ell^2+1))\bar f+4\bar
f^2+\frac{\ell^2}{4}.
 \eean
 In the notation (6.0.2), we have $$
a_1=16,\quad a_2=4(\ell^2+1),\quad a_3=0,\quad
c=-\frac{\ell^2}{4}. $$ Solving (6.0.3) for
$\alpha,\beta,\gamma,\delta$ leads to the canonical
form for Painlev\'e V, with $$ \alpha=-\beta
=\frac{(1+\ell )^{2}}{8},\quad \gamma =0,\quad\delta
=-8, $$ and according to Appendix 4,
 $f_{\ell}^{\prime\prime}(0)=2$, ending the proof of the first half
of Theorem 0.1.

Of course, from combinatorics (Proposition 1.4), we
have a much stronger statement:
 $$
 E_{O_{\pm}(\ell  +1)} e^{x \Tr M} =
\exp \left(
 \frac{x^2}{2}\pm\frac{x^{\ell+1}}{(\ell+1)!}
 +O(x^{\ell+2})\right),
 $$
and thus
  $$ f_{\ell}^{\pm}(x)=x\frac{d}{dx} \log
 E_{O(\ell+1)_{\pm}}e^{x\Tr M}dM=x^2 \pm
\frac{x^{\ell+1}}{\ell !}+O(x^{\ell+2}),~\mbox{near}~
x=0 .$$

It remains to show the uniqueness of the solution to
the initial value problem (3.2.2). Indeed,
substituting $f(x)=x^2+\sum_{i\geq3}a_ix^i$ into the
third order differential equation (3.2.2) for $f$
yields the recursive formula for the coefficients:
 $$%a_2=1,~~
 3(4-\ell^2)a_3=0$$
  \be
(i+1)(i^2-\ell^2)a_{i+1}-16(i-2)a_{i-1}+\sum_{2 \leq
m,n \leq i-1 \atop n+m=i+1} n
a_{n}(6m-4)a_m=0,~~\mbox{for}~ i\geq 3. \ee
 Therefore, if $\ell \geq 3$, we have inductively $a_3=...=a_{\ell}=0$, from (3.2.5). Setting $i=\ell$,
 in the equation above shows that the coefficient
 $a_{\ell+1}$ is free and can therefore be specified;
 it is specified by the combinatorics, namely
  $a_{\ell +1}=\pm ((\ell +1)!)^{-1}$.
 Once $a_{\ell+1}$ is fixed, all the subsequent $a_i$'s
 are determined by (3.2.5).

 \qed

\section{Painlev\'e equations for $U(n)$ integrals}

%\subsection{Toeplitz, Virasoro and Painlev\'e equations }

%\newpage

In this section, we show  items (ii) and (iii) of
Theorem 0.1:
\begin{proposition}
 \be
 g_n(x)=\frac{d}{dx} x \frac{d}{dx}\log
\int_{U(n)}e^{\sqrt{x}\Tr (M+\bar M)} dM
 \ee
 is the unique solution to the initial value problem
(Painlev\'e V equation):
 $$
  \mbox{\bf}
\left\{\begin{array}{l} \displaystyle{
g^{\prime\prime}_n-\frac{g^{\prime
2}_n}{2}\left(\frac{1}{g_n-1}+\frac{1}{g_n}
\right)+\frac{g^{\prime}_n}{x}-\frac{n^2}{2x^2}
 \frac{(g_n-1)}{g_n}
 +\frac{2}{x}g_n(g_n-1)
=0
}
    \\  \\
\displaystyle{\mbox{with}~~ g_{n}(x)=1-\frac{x^{n}}{(n
!)^2}+O(x^{n+1}) ,~\mbox{near}~ x=0 .}
\end{array}
 \right.
  $$

\end{proposition}

\begin{proposition}
\bea
  h_{n}(x)&=&\frac{E_{U(n)}\Tr M\det(I+M)^k e^{-x\Tr
\bar M} }{E_{U(n)}\det(I+M)^k e^{-x\Tr  \bar M} }
 \nonumber\\
 &=&\frac{1}{n+k}x\frac{d}{dx}\log E_{U(n)}
\det (I+M)^ke^{-x \Tr(I+\bar M)}dM
 \eea is the unique solution to the initial value Painlev\'e V equation, as
well:
 \be
\mbox{} \left\{\begin{array}{l} \displaystyle{
h^{\prime\prime\prime}-\frac{1}{2}\left(\frac{1}{h'}+
\frac{1}{h'+1}\right)h^{\prime\prime
2}+\frac{h''}{x}+\frac{2(n+k)}{x}h'(h'+1)}%\nonumber
 \\
\displaystyle{ -\frac{1}{{2x^2
h^{\prime}(h^{\prime}+1)}} \Bigl((x-n)h^{\prime}-h-n
\Bigr)
 \Bigl((2h +x+n)h^{\prime}+h+n \Bigr) =0
 }%\nonumber
 \\
    \\
\displaystyle{\mbox{with}~~
h:=h_{n}(x)=x\frac{k-n}{k+n}-\frac{x^{n +1
}}{(n+1)!}\left( k+n-1 \atop n \right)
  +O(x^{n+2}) ,~\mbox{near}~ x=0 .}
\end{array}
 \right. \nonumber\ee

\end{proposition}

\noindent\underline{\sl Proof of Propositions 4.1}:
The proofs of the two propositions are almost
identical,
 except in the end one specializes to a different locus.

  Throughout, we shall be using the
diagonal elements (2.1.12) of $L_1$ and
$hL^{\top}_2h^{-1}$:
 \be
 b_n=\frac{\pl}{\pl t_1}\log
 \frac{\tau_n}{\tau_{n-1}}=(L_1)_{n-1,n-1}
 ~~~\mbox{and}~~~
 b^*_n=-\frac{\pl}{\pl s_1}\log \frac{\tau_n}{\tau_{n-1}}
 =(hL_2^{\top}h^{-1})_{n-1,n-1}.
 \ee
 From (2.4.2), (2.1.17), (2.1.3), Theorem 2.5, and (5.0.4),
the integral below, which is also the determinant of a Toeplitz
matrix,
 \bea
\tau_n(t,s)&=&\int_{U(n)}e^{\sum_1^{\iy}\Tr
(t_iM^i-s_i\bar M^i)} dM\nonumber \\ &=&
 \det\left(\int_{S^1}z^{k-\ell}
e^{\sum_1^{\iy}(t_i z^i-s_iz^{-i})}
 \frac{dz}{2\pi i z}\right)_{0\leq k,\ell\leq n-1}
  \eea
satisfies the following three relations:

\noindent (i) {\bf Toeplitz}:
\begin{eqnarray}
{\cal T}(\tau)_{n}&=&\frac{\pl}{\pl t_{1}}\log
\frac{\tau_n}{\tau_{n-1}}
\frac{\pl}{\pl
s_{1}}\log\frac{\tau_n}{\tau_{n-1}}\nonumber\\
 && \hspace{0cm}+
 \left(1+\frac{\pl^{2}}{\pl s_{1}\pl t_{1}}
  \log\tau_n\right)\left(1+\frac{\pl^{2}}{\pl s_{1}\pl
t_{1}}\log\tau_n-\frac{\pl}{\pl s_{1}}
\left(\frac{\pl}{\pl t_{1}}\log\frac{\tau_{n}}
{\tau_{n-1}}\right)\right)
 \nonumber\\
 &=&-b_n b^*_n+\left( 1+\frac{\pl^{2}}{\pl s_{1}\pl t_{1}}\log\tau_n
   \right)  \left(1+ \frac{\pl^{2}}{\pl s_{1}\pl t_{1}}\log\tau_n
   -\frac{\pl}{\pl s_1}b_n  \right)=0, \nonumber\\
\end{eqnarray}

\noindent (ii) {\bf two-Toda}:
\begin{eqnarray}
\frac{\pl^2\log\tau_{n}}{\pl s_2\pl
t_1}&=&-2\frac{\pl}{\pl
s_1}\log\frac{\tau_{n}}{\tau_{n-1}}\frac{\pl^{2}}{\pl s_{1}\pl
t_{1}}\log\tau_{n}- \frac{\pl^{3}}{\pl s_{1}^2\pl
t_{1}}\log\tau_{n} \nonumber\\
 &=& 2b_n^* ~\frac{\pl^{2}}{\pl s_{1}\pl
t_{1}}\log\tau_{n}- \frac{\pl^{3}}{\pl s_{1}^2\pl
t_{1}}\log\tau_{n},
\end{eqnarray}

\noindent (iii) {\bf Virasoro}:
\begin{eqnarray}
{\cal V}_{-1}\tau_n&=&\left(\sum_{i\geq
1}(i+1)t_{i+1}\frac{\pl}{\pl t_{i}}-\sum_{i\geq
2}(i-1)s_{i-1}\frac{\pl}{\pl s_{i}}+n\left(t_1+\frac{\pl}{\pl
s_{1}}\right)\right)\tau_n=0\nonumber\\ {\cal
V}_{0}\tau_n&=&\sum_{i\geq 1}\left(it_{i}\frac{\pl}{\pl
t_{i}}-is_{i}\frac{\pl}{\pl s_{i}}\right)\tau_n=0\nonumber\\ {\cal
V}_{1}\tau_n&=&\left(-\sum_{i\geq 1}(i+1)s_{i+1}\frac{\pl}{\pl
s_{i}}+\sum_{i\geq 2}(i-1)t_{i-1}\frac{\pl}{\pl
t_{i}}+n\left(s_1+\frac{\pl}{\pl t_1}
\right)\right)\tau_n=0.\nonumber\\
\end{eqnarray}
Therefore we have
\begin{eqnarray}
0&=&\frac{1}{\tau_{n}}({\cal V}_{-1}+{\cal V}_{0})\tau_n\nonumber
\\
&=&\Bigl(\sum_{i\geq 1}((i+1)t_{i+1}+it_{i})\frac{\pl}{\pl t_i}-
\sum_{i\geq 2} ((i-1)s_{i-1}+is_{i})\frac{\pl}{\pl
s_i}\nonumber\\ && ~~~+(n-s_{1})\frac{\pl}{\pl
s_{1}}\Bigr)\log\tau_n+nt_1\nonumber\\ 0&=&\frac{1}{\tau_{n}}({\cal
V}_{0}+{\cal V}_{1})\tau_n\nonumber\\ &=&\Bigl(\sum_{i\geq
2}((i-1)t_{i-1}+it_{i})\frac{\pl}{\pl t_{i}}-\sum_{i\geq
1}((i+1)s_{i+1}+is_{i})
 \frac{\pl}{\pl s_{i}}\nonumber\\
 &&~~~~+
(n+t_{1})\frac{\pl}{\pl t_{1}}\Bigr)\log\tau_n+ns_{1}\nonumber\\
0&=&\frac{\pl}{\pl t_{1}}\left(\frac{{\cal
V}_{-1}\tau_{n}}{\tau_{n}}\right)\nonumber\\
 &=&\left(\sum_{i\geq
1}(i+1)t_{i+1}\frac{\pl^{2}}{\pl t_{1}\pl t_{i}}- \sum_{i\geq 2}
(i-1)s_{i-1}\frac{\pl^2}{\pl t_{1}\pl s_{i}}+n\frac{\pl^2}{\pl
t_1\pl s_1}\right)\log\tau_n+n\nonumber\\ 0&=&\frac{\pl}{\pl
t_{1}}\left(\frac{{\cal V}_{0}\tau_{n}}{\tau_{n}}\right)\nonumber\\
&=&\left(\sum_{i\geq 1}\left(it_{i}\frac{\pl^{2}}{\pl t_{1}\pl
t_{i}}-is_{i}
\frac{\pl^2}{\pl t_{1}\pl s_{i}}\right)+\frac{\pl}{\pl
t_{1}}\right)\log\tau_n  \nonumber\\
     0&=&\frac{\pl}{\pl s_{1}}\left(\frac{{\cal
V}_{1}\tau_{n}}{\tau_{n}}\right)\nonumber\\
 &=&\left(-\sum_{i\geq
1}(i+1)s_{i+1}\frac{\pl^{2}}{\pl s_{1}\pl s_{i}}+ \sum_{i\geq 2}
(i-1)t_{i-1}\frac{\pl^2}{\pl s_{1}\pl t_{i}}+
 n\frac{\pl^2}{\pl
s_1\pl t_1}\right)\log\tau_n+n\nonumber\\
 0&=&\frac{\pl}{\pl
s_{1}}\left(\frac{{\cal V}_{0}\tau_{n}}{\tau_{n}}\right)\nonumber\\
&=&\left(-\sum_{i\geq 1}\left(is_{i}\frac{\pl^{2}}
 {\pl s_{1}\pl
s_{i}}-it_{i}
\frac{\pl^2}{\pl s_{1}\pl t_{i}}\right)-
 \frac{\pl}{\pl s_{1}}\right)\log\tau_n\nonumber\\
\end{eqnarray}

\medbreak

%\newpage

%\noindent\underline{Case 1}.
For the sake of this proof, consider the
 $$\mbox{locus
${\cal L}=\{~\mbox{all}~t_{i}=s_{i}=0$,except
$t_{1},s_{1}\neq 0\}$.}$$

\noindent From (4.0.7), we have on $\LR$,
$$
\frac{{\cal V}_{0}\tau_{n}}{\tau_{n}}\Big|_{\LR}
=\left(t_{1}\frac{\pl}{\pl t_{1}}-s_{1}\frac{\pl}{\pl
s_{1}}\right)\log\tau_{n}\Big|_{\LR}=0,
$$
implying $\tau_{n}(t,s)\Big|_{\LR}$ is a function of
$x:=-t_{1}s_{1}$ only. Therefore we may write
$\tau_{n}\Big|_{\LR}=\tau_{n}(x)$, and so, along $\LR$, we have\\
$$
\frac{\pl}{\pl t_{1}}=-s_{1}\frac{\pl}{\pl x},\quad\frac{\pl}{\pl
s_{1}}=-t_{1}\frac{\pl}{\pl x},\quad\frac{\pl^2}{\pl t_{1}\pl
s_{1}}=-\frac{\pl}{\pl x}x\frac{\pl}{\pl x}.
$$
Setting
$$
%g_{n}(x)=\frac{\pl}{\pl
%x}\log\frac{\tau_{n}(t,s)}{\tau_{n-1}(t,s)}\Big|_{\LR}\mbox{\,\,and\,\,}
f_{n}(x)=
\frac{\pl}{\pl x}x\frac{\pl}{\pl x}\log\tau_{n}(x)=-
\frac{\pl^2}{\pl t_{1}\pl s_{1}}\log\tau_{n}(t,s)\Big|_{\LR},
$$
and using $x=-t_1s_1$, the two-Toda relation (4.0.6) takes on the form
\bea
s_1\frac{\pl^2\log\tau_n}{\pl s_2\pl t_1}\Big|_{\LR}
 &=&s_1\left(2b_n^* ~\frac{\pl^2}{\pl s_1\pl t_1}
  \log\tau_n-
\frac{\pl}{\pl s_1}\left(\frac{\pl^2\log\tau_n}
{\pl s_1\pl t_1}\right)\right)\nonumber\\
&=&x(2\frac{b_n^*}{t_1}f_n+f^{\prime}_n).
\eea
Setting this relation (4.0.9) into the Virasoro
 relations (4.0.7) and (4.0.8), we have
\begin{eqnarray}
0=\frac{{\cal V}_{0}\tau_{n}}{\tau_{n}}
-\frac{{\cal V}_{0}\tau_{n-1}}{\tau_{n-1}}\Big|_{\LR}
&=&\left(t_{1}\frac{\pl}{\pl t_{1}}-s_{1}\frac{\pl}{\pl
s_{1}}\right)\log\frac{\tau_{n}}{\tau_{n-1}}\Big|_{\LR}=
 t_1b_n+s_1b_n^* \nonumber\\ &&\\
0=  \frac{\pl}{\pl t_{1}}\frac{{\cal
V}_{-1}\tau_{n}}{\tau_{n}}\Big|_{\LR} &=&
\left(-s_{1}\frac{\pl^2}{\pl s_{2}\pl t_{1}}+n \frac{\pl^2}{\pl
t_{1}\pl s_{1}}\right)\log\tau_{n}\Big|_{\LR}+n \nonumber
\\ & &\nonumber\\
 &=& -x\left( 2\frac{b_n^*}{t_1} f_n(x)+f'_n(x)
   \right) +n (-f_n(x)+1).\nonumber\\
\end{eqnarray}
This is a system of two linear relations (4.0.10) and (4.0.11) in
$b_n$ and $b_n^*$, whose solution, together with its derivatives,
are given by: $$ \frac{b_n^*}{t_1}=-\frac{b_n}{s_1}=-
 \frac{n(f_n-1)+xf^{\prime}_n}{2xf_n},\quad
\frac{\pl b_n}{\pl s_1}=\frac{\pl }{\pl
x}x\frac{b_n}{s_1}=\frac{x(f_n
f^{\prime\prime}_n-f^{\prime
2}_n)+(f_n+n)f^{\prime}_n}{2f^2_n}.
 $$
  Setting $\pl^2 \log \tau_n /\pl s_1\pl t_1 =-f_n$
  into the Toeplitz relation (4.0.5) yields
  $$
b_nb_n^* =(1-f_n)\left(1-f_n-\frac{\pl}{\pl
s_1}b_n\right),
 $$
 which, using the expressions above for $b_n,~b_n^{*}$
 and $\pl b_n/\pl s_1$, yields the differential
 equation:
 \be
f^{\prime\prime}_n-\frac{1}{2}f^{\prime
2}_n\left(\frac{1}{f_n-1}+\frac{1}{f_n}
\right)+\frac{1}{x}f^{\prime}_n+\frac{n^2(-f_n+1)}{2x^2f_n}
-\frac{2}{x}f_n(-f_n+1) =0.
 \ee
  Note, along the locus $\LR$, we may set
  $t_{1}=\sqrt{x}$ and
$s_{1}=-\sqrt{x}$, since it respects $t_1s_{1}=-x$.
Thus,
 $$
f_{n}(x)= \frac{d}{d x}x\frac{d}{d
x}\log\tau_{n}(x)%=- \frac{\pl^2}{\pl t_{1}\pl s_{1}}\log\tau_{n}(t,s)\Big|_{\LR},
, $$
 with
 $$
\left.\tau_{n}(t,s)\right|_{\LR}=\left.\int_{U(n)}
e^{\Tr(t_1M-s_1\bar M)}dM\right|_{\LR}=
\int_{U(n)}e^{\sqrt{x}\,\Tr(M+\bar M)}dM ,
 $$
  satisfies
(4.0.12). The behavior of $f_n(x)$ near $x=0$ is given
by Proposition 1.5 and the above formula, with the
uniqueness established as in the orthogonal case,
proving Proposition 4.1. \qed

 \remark Setting $$ f_n(x)=\frac{w(x)}{w(x)-1} $$ leads to
standard Painlev\'e V, with $\alpha=\delta=0,~ \beta=-n^2/2,
\gamma=-2$.

\bigbreak

%\newpage

\noindent\underline{\sl Proof of Propositions 4.2}:
%\noindent\underline{Case 2}.
For fixed $k\in \BR,
~k\neq 0$, consider the locus $$\LR=\{\mbox{all}~
it_i=-k(-1)^i ~\mbox{ and}~ s_i=0 ,~\mbox{except}~
s_1=x\}.$$ Then, setting
 \be
f_n(x)=\frac{\pl}{\pl t_1}\log\tau_n\Big|_{\LR}, \ee
the Toda relations (4.0.6) become:
\begin{eqnarray}
-x\frac{\pl^2}{\pl s_2\pl t_1}\log\tau_{n}\Big|_{\LR}
 &=& -2xb_n^* ~\frac{\pl^{2}}{\pl s_{1}\pl
t_{1}}\log\tau_{n}+ x\frac{\pl^{3}}{\pl s_{1}^2\pl
t_{1}}\log\tau_{n}\nonumber\\
  &=&- 2xb_n^*~f_n'+xf_n''.
\end{eqnarray}
The Virasoro relations (4.0.8) become, using (4.0.3) and the
locus,
  \bea 0&=&\left(\frac{({\cal V}_0+{\cal
V}_1)\tau_n}{\tau_n}- \frac{({\cal V}_0+{\cal
V}_1)\tau_{n-1}}{\tau_{n-1}}\right)\Big|_{\LR}\nonumber\\
&=&\left(-x\frac{\pl}{\pl s_1}+(n+k)\frac{\pl}{\pl
t_1}\right)\log\tau_n+nx\nonumber\\ &
&\hspace{1cm}-\left(-x\frac{\pl}{\pl
s_1}+(n-1+k)\frac{\pl}{\pl t_1}\right)
\log\tau_{n-1}-(n-1)x\nonumber\\
&=&\left(-x\frac{\pl}{\pl s_1}+(k+n-1)\frac{\pl}{\pl
t_1}\right)\log\frac{\tau_n}{\tau_{n-1}}+\frac{\pl}{\pl
t_1}\log\tau_n+x \nonumber\\ &=&-x\frac{\pl}{\pl
s_1}\log\frac{\tau_n}{\tau_{n-1}}+(k+n-1)\frac{\pl}{\pl
t_1}\log\frac{\tau_n}{\tau_{n-1}}+f_n+x\nonumber \\
 &=& xb_n^*+(k+n-1)b_n+f_n+x  ,
\eea and, using (4.0.8) and (4.0.14) \bea
0&=&\frac{\pl}{\pl t_1}\frac{({\cal V}_{-1}+{\cal
V}_0)\tau_n}{\tau_n}\Big|_{\LR}\nonumber\\
&=&\left(\frac{\pl}{\pl t_1}-x\frac{\pl^2}{\pl t_1\pl
s_2}+(n-x)\frac{\pl^2}{\pl s_1\pl
t_1}\right)\log\tau_n+n=0\nonumber\\
&=&f_n+(n-x)f_n'+n-2xb_n^* f_n'+xf_n''. \nonumber\\
\eea
 So,
as before, we have a linear system in $b_n$ and
$b_n^*$ whose solution is
 \bean
 b_n &=&-\frac{x
f_n^{\prime\prime}+f_n^{\prime}(2f_n + x+n)+f_n+n}{2
f_n^{\prime}(n+k-1)} \\
 b^*_n &=&\frac{x f_n^{\prime\prime}- f_n^{\prime}(x-n)+f_n+n}
  {2 x f_n^{\prime}}.
 \eean
Substituting this solution into the Toeplitz relation
(4.0.5): $$ b_nb_n^*=(1+f_n')(1+f_n'-\frac{\pl}{\pl
x}b_n)
 $$
  yields
\bean
&&f_n^{\prime\prime\prime}-\frac{1}{2}\left(\frac{1}{f_n'}+
\frac{1}{f_n'+1}\right)f_n^{\prime\prime
2}+\frac{f_n''}{x}+\frac{2(n+k)}{x}f_n'(f_n'+1)\\ &&
-\frac{1}{{2x^2f_n^{\prime}(f_n^{\prime}+1)}}
\Bigl((x-n)f_n^{\prime}-f_n-n \Bigr)
 \Bigl((2f_n +x+n)f_n^{\prime}+f_n+n \Bigr) =0.
\eean
 It remains to compute $f_n(x)$ as in (4.0.13). Note that
  \bean
 \tau_n(x):=\left.\tau_n(t,s)\right|_{\LR}&=&
 \left. \int_{U(n)}e^{\Tr\sum_{1}^{\iy}(t_iM^i-s_i\bar
M^i)}dM\right|_{\LR}
 \\
 &=&\int_{U(n)}
\left(e^{-\Tr\sum_{1}^{\iy}\frac{(-M)^i}{i}} \right)^k
e^{-x\Tr\bar M}dM
 \\
 &=&\int_{U(n)} \det(I+M)^ke^{-x\Tr\bar M}dM.
 \eean
 Therefore \bea f_n(x)&=&\left.\frac{\pl}{\pl
t_{1}}\log\tau_n\right|_{\LR}\nonumber\\
&=&\left.\frac{\int \Tr
M\,e^{\Tr\sum_{1}^{\iy}(t_iM^i-s_i\bar M^i)}dM} {\int
e^{\Tr\sum_{1}^{\iy}(t_iM^i-s_i\bar
M^i)}dM}\right|_{\LR} \nonumber\\ &=&\frac{\int \Tr
M\det(I+M)^ke^{-x \Tr\bar M}dM} {\int \det(I+M)^k
e^{-x \Tr\bar M} dM}\nonumber\\ &\stackrel{*}{=}&
\frac{1}{n+k} x \frac{d}{dx}\log \tau_n e^{-nx}. \eea
% This proves Proposition 4.2; except for this last
%equality, which will be shown in the lemma below, and so the proof of
%Theorem 0.1 is finished.
This last equality $\stackrel{*}{=}$ will be shown
later in Lemma 4.3. To conclude the proof of
Proposition 4.2, observe from (4.0.17) and Proposition
1.5, that
 \bean
 f_n(x)&=& \frac{1}{n+k} \left(x \frac{d}{dx}\log
\int_{U(n)}\det(I+M)^k e^{-x \Tr \bar M} d M
-nx\right)\\
  &=&
 x\frac{k-n}{k+n}- \frac{x^{n+1}}{(n+1)!}
  \left( k+n-1 \atop n
\right)
  +O(x^{n+2}),
 \eean
 concluding the proof of Proposition 4.2.\qed

\noindent\underline{\sl Proof of Theorem 0.1}: Upon
integrating the expressions (4.0.1) and (4.0.2) and
exponentiating, one finds the expressions (ii) and
(iii) of Theorem 0.1, upon using respectively the
initial conditions (0.0.3) and the first identity of
Lemma 4.3. \qed

Recall equality (4.0.17) ($\stackrel{*}{=}$) still
needed proof:

\begin{lemma}
 $$ \tau_n(0)=\int_{U(n)}\det(I+M)^k dM=1 ~~\mbox{,}~~
 \frac{\pl \tau_n}{\pl t_1}(0)= \int_{U(n)}\Tr M \det(I+M)^k dM=0 $$
and
 \bea %\hspace{-.2cm}
 \lefteqn{\frac{ \int
\Tr M\det(I+M)^ke^{-x \Tr\bar M}dM} { \int \det(I+M)^k
e^{-x\Tr \bar M} dM}}
 \nonumber\\
  &=&\frac{1}{n+k} x
\frac{d}{dx}\log
  e^{-nx}\int_{U(n)} \det(I+M)^ke^{-x \Tr \bar M}dM
  \nonumber\\
 &=&\frac{-x}{n+k}\left(\frac{\int\Tr \bar
M\det(I+M)^ke^{-x \Tr \bar M}dM} {\int\det(I+M)^ke^{-x
\Tr\bar M}dM}+n\right). \eea

\end{lemma}

\proof  Recall from (4.0.17) that $f(x)$ is the left
hand side of (4.0.18). At first we show, using the
Toeplitz matrix representation (2.3.11) in the third
identity, that $f(0)=0$; indeed,
 \bea
  \tau_n(0)~f(0)&=&\int_{U(n)} \Tr M\det(I+M)^kdM
   \nonumber\\
&=&\frac{d}{d\varepsilon}\left.\int\det(I+M)^k\det(I+\varepsilon
M)dM\right|_{\varepsilon =0}
 \nonumber\\
&=&\frac{d}{d\varepsilon}\det\left.\left(\int_{S^1}z^{\ell
-m}(1+z)^k(1+\varepsilon z)
\frac{dz}{2\pi\,iz}\right)_{0\leq\ell,m\leq
n-1}\right|_{\varepsilon =0}
 \nonumber\\
&\stackrel{*}{=}&\frac{d}{d\varepsilon}\det\left.\left(
\begin{array}{cccc}
1& & & \\
 &1&*& \\
 & &\ddots& \\
 O& & &1
\end{array}\right)
\right|_{\varepsilon =0}
 \nonumber\\
 &=&\frac{d}{d\varepsilon}(1)=0. \eea
  The equality
$\stackrel{*}{=}$ is due to the fact that

 $$
\begin{array}{lll}
\int_{S^1} z^{\ell -m}(1+z)^k(1+\varepsilon
z)\displaystyle{\frac{dz}{2\pi\,iz}}&=0&\mbox{for $\ell
-m\geq 1$}\\
 \\
 &=1&\mbox{for $\ell =m$.}
\end{array}
$$
 The same, but even simpler argument shows
$\tau_n(0)=1$, by replacing $1+\varepsilon z$ by $1$
in (4.0.19). From (4.0.8), it also follows that
 \bean
0&=&\frac{\pl}{\pl s_{1}}
\left.\frac{(\VR_0+\VR_1)\tau_n}{\tau_n}\right|_{\LR}
 \\
&=&\left.\left((n+k)\frac{\pl^2}{\pl s_{1}\pl
t_{1}}-s_{1}\frac{\pl^2}{\pl s_{1}^2}-\frac{\pl}{\pl
s_{1}}\right)\log\tau_n \right|_{\LR} +n
 \\
&=&(n+k)\left.\frac{\pl f}{\pl x}-\frac{\pl}{\pl x}
x\frac{\pl}{\pl x}\log\tau_n \right|_{\LR} +n
 \\
&=&(n+k)\left.\frac{\pl f}{\pl x}-\frac{\pl}{\pl
x}x\frac{\pl}{\pl x} \log\tau_n e^{-nx}\right|_{\LR}.
 \eean
 Integrating this expression from $0$ to $x$
yields $$ (n+k)(f(x)-f(0))=x\frac{\pl}{\pl
x}\log\tau_n e^{-nx};
  $$
  the fact that $f(0)=0$ establishes the first identity of
  (4.0.18). The second identity of (4.0.18)
  follows from
   \bean
f(x) &=&\frac{1}{n+k} x \frac{d}{dx}\log
  e^{-nx}\int_{U(n)} \det(I+M)^ke^{-x\Tr \bar M}dM\\
   &=&\frac{-x}{n+k}\left(\frac{\int \Tr \bar
M\det(I+M)^ke^{-x \Tr \bar M}dM} {\int\det(I+M)^ke^{-x
\Tr \bar M}dM}+n\right),
 \eean
 ending the proof of Lemma 4.3.\qed

\section{Appendix 1: Virasoro algebras}

In \cite{AvM3}, we defined a Heisenberg and Virasoro algebra of
vector operators $^{\beta}\BJ^{(i)}_{k}$, depending on a parameter
$\beta
>0$:
$$
\left( ^{\beta}\BJ_k^{(1)}\right)_n=\,\,^{\beta}J_k^{(1)}+nJ_k^{(0)}
 ~\mbox{and}~ \left(\BJ_k^{(0)}\right)_n=nJ_k^{(0)}= n\dt_{0k}
,$$
 and
\vspace{0.3cm}

\noindent$\displaystyle{{}^{\beta}\BJ_k^{(2)}(\beta)}$
\bea
&=&\frac{\beta}{2}\sum_{i+j=k}
:{}^{\beta}\BJ_i^{(1)}{}^{\beta}\BJ_j^{(1)}:
+\left(1-\frac{\beta}{2}\right)\left((k+1) \,\, ^{\beta}\BJ_k^{(1)}
-k~\BJ_k^{(0)}\right)\nonumber\\
&=&\left(\frac{\beta}{2}. ~{}^{\beta}J_k^{(2)} +
\left(n\beta +(k+1)(1-\frac{\beta}{2})\right).
 ~{}^{\beta}J_k^{(1)}
+ \frac{n((n-1)\beta+2)}{2} J_k^{(0)}\right)_{n\in \BZ} .
 \nonumber\\
 \eea
The ${}^{\beta}\BJ_k^{(2)}$'s satisfy the commutation relations:
(see \cite{AvM3})
\begin{eqnarray}
\left[~ {}^{\beta}\BJ_k^{(1)},{}^{\beta}\BJ_{\ell}^{(1)}
  \right] &=&\frac{k}{\beta}\delta_{k,-\ell}\nonumber\\
\left[~{}^{\beta}\BJ_k^{(2)},~{}^{\beta}\BJ_{\ell}^{(1)} \right]
&=&-\ell ~ ~{}^{\beta}\BJ_{k+\ell}^{(1)}
+k(k+1)\left(\frac{1}{\beta}-
\frac{1}{2}\right)\delta_{k,-\ell}\nonumber \\
\left[~{}^{\beta}\BJ_k^{(2)},~{}^{\beta}\BJ_{\ell}^{(2)}
 \right]
&=&(k-\ell)~{}^{\beta}\BJ_{k+\ell}^{(2)} +c\left( \frac{k^3-k}{12}
\right)\dt_{k,-\ell}~,
\end{eqnarray}
with central charge
$$
c=1-6\left(\left(\frac{\beta}{2}\right)^{1/2}
-\left(\frac{\beta}{2}\right)^{-1/2}  \right)^2.
$$
In the expressions above,
\bea
 {}^{\beta}J_k^{(1)}&=&\frac{\pl}{\pl
t_k}~\mbox{for}~~k>0
 \nonumber\\
&=&\frac{1}{\beta}(-k)t_{-k}~\mbox{for}~~k<0\nonumber\\
 &=&0~\mbox{for}~~k=0\nonumber\\  \nonumber\\
 ^{\beta}J^{(2)}_{k}&=&\sum_{i+j=k}\frac{\pl^2}{\pl
 t_{i}\pl t_{j}}+\frac{2}{\beta}\sum_{-i+j=k}it_{i}\frac{\pl}{\pl
 t_{j}}+\frac{1}{\beta^2}\sum_{-i-j=k}it_{i}jt_{j}
 \eea
In particular, for $\beta = 1$ and $2$, the
 ${}^{\beta}\BJ_k^{(2)}$ take on the form:
\bea
{}^{\beta}\BJ_k^{(2)}(t)\Big|_{\beta=1}
 %=\left(J_{k,n}^{(2)}\right)_{n\in\BZ}
 &=&\left.\frac{1}{2}\left(~{}^{\beta}J_k^{(2)}
+(2n+k+1)\,\, ^{\beta}J_k^{(1)} +n(n+1)J_k^{(0)}
 \right)_{n\in\BZ}\right|_{\beta=1},\nonumber\\   &&\\
 {}^{\beta}\BJ_k^{(2)}(t)\Big|_{\beta=2}
 %=\left(J_{k,n}^{(2)}\right)_{n\in\BZ}
 &=&\left.\left(~{}^{\beta}J_k^{(2)}
+2n\,\, ^{\beta}J_k^{(1)} +n^2 J_k^{(0)}
 \right)_{n\in\BZ}\right|_{\beta=2}
\eea

\section{Appendix 2: Chazy classes}

Given arbitrary polynomials $P(z), Q(z), R(z)$ of degree $3,2,1$
respectively, Cosgrove \cite{C}, (A.3), shows that the following
third order equation
\be
f^{\prime \prime\prime}+\frac{P'}{P}f^{
\prime\prime}+\frac{6}{P}f^{\prime
2}-\frac{4P'}{P^2}ff' +\frac{P^{\prime\prime}}{P^2}
f^2
+\frac{4Q}{P^2}f'-\frac{2Q'}{P^2}f+\frac{2R}{P^2}=0
\ee has a first integral, which is second order in $f$
and quadratic in $f^{\prime\prime}$, \bea &f^{
\prime\prime 2}& +\frac{4}{P^2}
 \left( (Pf^{\prime 2}+Q f^{\prime}+R)f^{\prime}
 - (P' f^{\prime 2}+\frac{}{}Q' f^{\prime}+R')f^{}
  \right. \nonumber\\ && \hspace{2cm} \left.
   +\frac{1}{2}(P^{\prime\prime}f^{\prime
}+Q^{\prime\prime} )f^2
 -\frac{1}{6} P^{\prime\prime\prime}f^3 +c\right)=0;
\eea $c$ is the integration constant. This is a master Painlev\'e
equation, containing the 6 Painlev\'e equations. When the
polynomials $P$, $Q$ and $R$ have the following form $$ P=x,\quad
Q=-\frac{a_{1}}{4}x^2,\quad R=-\frac{1}{4}(a _{2} x+a_{3}), $$
then equation (6.0.2) can be reduced to the Painlev\'e V equation
\cite{CS}, p.70:
 \be
w^{\prime\prime}=\left(\frac{1}{2w}+\frac{1}{w-1}\right)w^{\prime
2}-\frac{1}{x}w^{\prime}+\frac{(w-1)^2}{x^2}\left(\alpha
w+\frac{\beta}{w}\right)+\frac{\gamma w}{x}+\frac{\delta
w(w+1)}{w-1}, \ee
 with
  \bean a_{1} &=&-2\delta\\ a_{2}
&=&\frac{1}{4}\gamma^2+2\beta\delta -\delta(1-\sqrt{2\alpha})^2\\
a_{3} &=&\beta\gamma+\frac{1}{2}\gamma(1-\sqrt{2\alpha})^2\\ c
&=&-\frac{1}{32}\gamma^2((1-\sqrt{2\alpha})^2-2\beta)+
\frac{1}{32}\delta((1-\sqrt{2 \alpha})^2+2\beta)^2. \eean
%\newpage

\section{Appendix 3: The volume of the orthogonal and
 symplectic groups}

Selberg's integral (see Mehta \cite{M}, p 340), renormalized over
$[-1,1]$,

\bean
&&\int_{[-1,1]^n}
\Delta_{n}
(x)^{2\gamma}\prod^n_{j=1}(1-x_j)^{\alpha}(1+x_j)^{\beta}dx_j
\\&&  ~~~~~~~~~~~~~=2^{n(\alpha+\beta+\gamma(n-1)+1)}  \prod^{n-1}_{j=0}
\frac{\Gamma(\alpha+j\gamma+1)\Gamma(\beta+j\gamma+1)
\Gamma(\gamma+j\gamma+1)}
{\Gamma(\gamma+1)\Gamma(\alpha+\beta-\gamma+\gamma(n+j)+2)}
\\&&  ~~~~~~~~~~~~~=2^{n(n+\alpha+\beta)}
 \prod^n_{j=1}\frac{j!~ \Gamma(j+\alpha)\Gamma(j+\beta)}
 {\Gamma(n+j+\alpha+\beta)},~~\mbox{upon setting}~\gamma=1,
\eean leads to the value of $c_{2n}^{\pm}$ and
$c_{2n-1}^{\pm}$ in Theorem 1.1:
\begin{itemize}
  \item  $\alpha=-\beta=\pm \frac{1}{2} ~~~ \longrightarrow
 ~~ \int_{O(2n+1)_{\pm}}dM=2^{n^2}\displaystyle{\prod_{j=1}^n
   \frac{j!(j-1/2)\Gamma^2(j-1/2)}{(n+j-1)!}}$
  \item $\alpha=\beta=- \frac{1}{2} ~~~ \longrightarrow
 ~~ \int_{O(2n)_{+}}dM=\displaystyle{2^{n(n-1)}\prod_{j=1}^n
   \frac{j! \Gamma^2(j-1/2)}{(n+j-2)!}}$
  \item $\alpha=\beta= \frac{1}{2},~n\mapsto n-1 ~~~ \longrightarrow
 ~~ \int_{O(2n)_{-}}dM=\displaystyle{2^{n(n-1)}\prod_{j=1}^{n-1}
   \frac{j! \Gamma^2(j+1/2)}{(n+j-1)!}}$.
\end{itemize}

%\newpage

\section{Appendix 4: Direct evaluation of integrals over the
orthogonal group and their derivatives at $x=0$}

Refering to Theorem 3.1, formulae (3.1.3) and (3.1.4),
we evaluate $d/dx\log \tau_n(x)$ and $d^2/dx^2\log
\tau_n(x)$ directly from the integral representation,
not using the combinatorial interpretation of the
integrals.  To do this, we need the Aomoto extension
\cite{A} (see Mehta \cite{M}, p. 340) of Selberg's
integral:\footnote{where $Re\,\alpha$, $Re\,\beta
>-1$, $Re\,\gamma >
-\min\displaystyle{\left(\frac{1}{n}, \frac{Re\,\alpha
+1}{n-1}, \frac{Re\,\beta +1}{n-1}\right)}$}

\bea
 \la x_1...x_m\ra&:=&\displaystyle{\frac{\int^1_{0}\ldots \int^1_{0}x_1\ldots
x_m\left|\Delta(x)\right|^{2\gamma}\prod^n_{j=1}x^{\alpha}_{j}
(1-x_{j})^{\beta}dx_1...dx_n}{\int^1_{0}\ldots
\int^1_{0}\left|\Delta(x)\right|^{2\gamma}
\prod^n_{j=1}x^{\alpha}_{j}
(1-x_{j})^{\beta}dx_1...dx_n}}\nonumber\\
 &=&\prod^m_{j=1}\frac{\alpha
+1+(n-j)\gamma}{\alpha+\beta+2+(2n-j-1)\gamma}. \eea In
particular, setting $\gamma=1$, formula (8.0.1) implies
 \be \la
x_1\ra={{n+\alpha}\over{2\,n+\beta+\alpha}}~~ \mbox{and}~~ \la x_1 x_2\ra=
{{\left(n+\alpha-1\right)\,\left(n+\alpha\right)}\over{\left(2\,n+\beta+\alpha-1
 \right)\,\left(2\,n+\beta+\alpha\right)}}, \ee
  and from the identity ( see \cite{M}, p. 349 )
 $$
 \left(2\,n+\beta+\alpha+1\right)\la x_1^2\ra   = \left(2\,n+\alpha\right)\la
 x_1\ra-(n-1)\la x_1x_2\ra,
 $$
 we derive
\be
  \la
x_1^2\ra={{\left(n+\alpha\right)\,\left(3\,n^{2}+2\,\beta\,n+3\,\alpha\,n+\alpha
\,\beta+\alpha^{2}-1%
 \right)}\over{\left(2\,n+\beta+\alpha-1\right)\,\left(2\,n+\beta+\alpha\right)\
\left(%
 2\,n+\beta+\alpha+1\right)}}. \ee

 We now consider the following ratio of integrals:
 (remember $
\rho_{\alpha\beta}(z):=(1-z)^{\alpha}(1+z)^{\beta} $)
 \be
 \la y_1...y_m \ra_{[-1,1]}:=
\displaystyle{\frac{\int_{[-1,1]^n} y_1...y_m \Delta_{n} (y)^{2}
\prod^n_{k=1}\rho_{(\alpha,\beta)}(y_k)dy_k}{\int_{[-1,1]^n}
\Delta_{n} (y)^{2}\prod^n_{k=1}\rho_{(\alpha,\beta)}(y_k)dy_k}}.
 \ee
The relationship between the two integrals (8.0.1) and (8.0.4) is
obtained by setting
 $x_{j}=\frac{1-y_{j}}{2}$; so, we have
$$\la x_1 \ra=\frac{1}{2}(1-\la y_1\ra ),~  \la x_1x_2 \ra=
 \frac{1}{4}(1-2\la y_1\ra  + \la y_1y_2\ra ),~      \la
x_1^2 \ra=  \frac{1}{4}(1-2\la y_1\ra  + \la y_1^2\ra ).
 $$
So, these relations, upon using (8.0.2) and (8.0.3) and upon
setting $\alpha=\frac{a+b}{2},~\beta=\frac{a-b}{2}$, yield

\be \la y_1 \ra= \frac{-b}{a+2n} ~, ~~\la y_1 y_2\ra=
{{b^{2}-a-2\,n}\over{\left(a+2\,n-1\right)\,\left(a+2\,n\right)}}
\ee

  $$ \la y_1^2 \ra=
{{b^{2}(a+n)+n(a+2n)^2-(a+2n)}\over{%
 \left(a+2\,n-1\right)\,\left(a+2\,n\right)\,\left(a+2\,n+1\right)}}.
  $$
Hence, setting $$I_n^{(\al,\beta)}(x):=\int_{[-1,1]^n}\Dt_n(z)^2
\prod_1^n e^{2xz_k}\rho_{\alpha,\beta}(z_k)dz_k~, $$ we compute
for future use:
 \bea
 \gamma (n)&:=& 2\left. \frac{I_n^{''}}{I_n}\right|_{x=0}
  \nonumber\\
  &=& 8\la (\sum_1^n y_i)^2\ra  \nonumber\\
  &=& 8~(n\la y_1^2\ra+n(n-1) \la y_1y_2\ra)  \nonumber \\
  && \nonumber\\
  &=& 8n~(\la y_1^2\ra+(n-1) \la y_1y_2\ra)  \nonumber \\
  && \nonumber \\
  &=&  8n~
  {{ (a+2n)(b^2n+a+n)-b^2}
   \over{
 \left(a+2\,n-1\right)\,\left(a+2\,n\right)\,
 \left(a+2\,n+1\right)}}  .
 \eea

 Note this formula applies to a general
 $I_n^{(\al,\beta)}(x)$, where a combinatorial
 interpretation is absent. These considerations will
 now be applied to the orthogonal case. Indeed,
considering the special values of $\alpha$ and $\beta$
and thus for $a$ and $b$, we evaluate:
\begin{itemize}
  \item $a=-1,~~~~b=0~~: ~~~\gamma(n)=2$
  \item $a=1,~~~~~~b=0~~~: ~~~\gamma(n)=2$
  \item  $a=0,~~~~~b=1~~~~: ~~~\gamma(n)=4$
  \item $a=0,~~~~~b=-1~~: ~~~\gamma(n)=4$.
\end{itemize}
It is easily seen that
 $$
 \left( x\frac{d}{dx}\log \int e^{x \Tr M}
 dM\right)^{\prime\prime}=2\frac{(\int e^{x \Tr M}
 dM)^{\prime\prime}}{\int e^{x \Tr M}
 dM} - 2\left(\frac{(\int e^{x \Tr M} dM)^{\prime}}
 {\int e^{x \Tr M} dM}\right)^2
+O(x)$$ and so, using (3.2.4) and the fact that the
volume $\int dM$ does not vanish,
 $$
\left. \left( x\frac{d}{dx}\log \int e^{x \Tr M}
 dM\right)^{\prime\prime}\right|_{x=0}=
  2\left.\frac{(\int e^{x \Tr M}
 dM)^{\prime\prime}}{\int e^{x \Tr M}
 dM} \right|_{x=0}
 $$
Using $f_{\ell}'(0)=0$ to evaluate $I_n'(0)$ below and
using (1.0.2)
and (8.0.6), we now verify in each of the cases:
\begin{eqnarray*}
f^{\prime\prime}_{2n-1}(0)&=&
\left.\left(x\frac{d}{dx} \log\int_{O(2n)_+}e^{x \Tr
M}dM\right)^{''}\right|_{x=0}  \\ &=& 2\left.
\frac{I_n^{''(-\frac{1}{2},-\frac{1}{2})}}{I_n^{(-\frac{1}{2},-\frac{1}{2})}}
\right|_{x=0}=\left.\gamma(n)\right|_{a=-1,b=0}=2\\
f^{\prime\prime}_{2n-1}(0)&=&
\left.\left(x\frac{d}{dx}\log\int_{O(2n)_-}e^{x \Tr
M}dM\right)^{\prime\prime}\right|_{x=0}\\
 &=&2
\frac{I_{n-1}^{\prime\prime(\frac{1}{2},\frac{1}{2})}}
{I_{n-1}^{(\frac{1}{2},\frac{1}{2})}}=
\left.\gamma(n-1) \right|_{a=1,b=0}=2\\
f^{\prime\prime}_{2n}(0)&=&
\left.\left(x\frac{d}{dx}\log \int_{O(2n+1)_+}e^{x \Tr
M}dM\right)^{''}\right|_{x=0}\\
&=&2\left.\frac{(e^xI_n^{(\frac{1}{2},-\frac{1}{2})
})''}{
e^xI_n^{(\frac{1}{2},-\frac{1}{2})}}\right|_{x=0}\\
&=&2\left.\left(\frac{I_{n}^{\prime\prime}+2I^{'}_{n}+I_{n}}{I_{n}}\right)\right|_{x=0}\\
&=&2\left(\frac{I_{n}^{\prime\prime}}{I_{n}}-1\right),\mbox{
using
$\left.\left(e^xI_{n}(x)\right)^{'}\right|_{x=0}$}
=I_{n}'(0)+I_{n}(0)=0\\
&=&-\left.2+\gamma(n)\right|_{a=0,b=+1}=2\\
f^{\prime\prime}_{2n}(0)&=& \left.\left(x\frac{d}{dx}
 \log\int_{O(2n+1)_-}e^{x \Tr
M}\right)^{\prime\prime}\right|_{x=0}=
2\left.\frac{(e^{-x} I_n^{(-\frac{1}{2},\frac{1}{2})
})''}
{e^{-x}I_n^{(-\frac{1}{2},\frac{1}{2})}}\right|_{x=0}\\
&=&2 \left.\frac{I_{n}^{''}-2I^{'}_{n}
+I_{n}}{I_{n}}\right|_{x=0}, \mbox{ using
$\left.\left(e^{-x}I_{n}(x)\right)^{'}\right|_{x=0}$}=0\\
&=&2\left.\left(\frac{I_{n}^{''}}{I_{n}}-1\right)=-2+\gamma(n)\right|_{a=0,b=-1}
=2.
 \end{eqnarray*}

\end{document}